\documentclass[12pt]{amsart}
\usepackage{amsmath,amscd,amssymb,amsfonts,graphics}
\usepackage{xcolor,hyperref}
\hypersetup{colorlinks=true,citecolor=black,linkcolor=black}
\newcommand{\hl}{\hyperlink}
\newcommand{\htt}{\hypertarget}
\newcommand{\h}{\hbox}
\newcommand{\q}{\quad}
\newcommand{\nin}{\noindent}
\newcommand{\bs}{\par\bigskip}
\newcommand{\ms}{\par\medskip}
\newcommand{\sk}{\par\smallskip}
\newcommand{\bsn}{\par\bigskip\noindent}
\newcommand{\msn}{\par\medskip\noindent}
\newcommand{\skn}{\par\smallskip\noindent}
\newcommand{\ges}{\geqslant}
\newcommand{\les}{\leqslant}
\newcommand{\1}{\hskip1pt}
\newcommand{\mcap}{\hbox{$\bigcap$}}
\newcommand{\mcup}{\hbox{$\bigcup$}}
\newcommand{\msqcup}{\hbox{$\bigsqcup$}}
\newcommand{\msum}{\hbox{$\sum$}}
\newcommand{\mopl}{\hbox{$\bigoplus$}}
\newcommand{\mprod}{\hbox{$\prod$}}
\newcommand{\A}{{\mathcal A}}

\newcommand{\Lc}{{\mathcal L}}

\newcommand{\OO}{{\mathcal O}}
\newcommand{\Sc}{{\mathcal S}}
\newcommand{\X}{{\mathcal X}}
\newcommand{\Y}{{\mathcal Y}}
\newcommand{\Zc}{{\mathcal Z}}
\newcommand{\Ab}{{\mathbb A}}
\newcommand{\PP}{{\mathbb P}}
\newcommand{\Q}{{\mathbb Q}}
\newcommand{\C}{{\mathbb C}}

\newcommand{\Lb}{{\mathbb L}}
\newcommand{\N}{{\mathbb N}}
\newcommand{\R}{{\mathbb R}}
\newcommand{\RR}{{\mathbf R}}
\newcommand{\Z}{{\mathbb Z}}
\newcommand{\ob}{{\mathbf 1}}
\newcommand{\Gr}{{\rm Gr}}
\newcommand{\al}{\alpha}
\newcommand{\be}{\beta}
\newcommand{\ga}{\gamma}
\newcommand{\Ga}{\Gamma}
\newcommand{\de}{\delta}
\newcommand{\De}{\Delta}
\newcommand{\la}{\lambda}
\newcommand{\si}{\sigma}
\newcommand{\Si}{\Sigma}
\newcommand{\ep}{\varepsilon}
\newcommand{\qh}{\widehat{q}}

\newcommand{\gat}{\widetilde{\gamma}}
\newcommand{\Ht}{{}\,\widetilde{\!H}{}}
\newcommand{\Pt}{{}\,\widetilde{\!P}{}}
\newcommand{\Xt}{{}\,\widetilde{\!\mathcal X}{}}
\newcommand{\Yt}{{}\,\widetilde{\!\mathcal Y}{}}
\newcommand{\Zt}{{}\,\widetilde{\!\mathcal Z}{}}

\newcommand{\ttt}{{}\,\widetilde{\!t}{}}
\newcommand{\kc}{\check{k}}
\newcommand{\Ah}{{}\,\widehat{\!A}{}}
\newcommand{\jh}{{}\,\widehat{\!j}{}}
\newcommand{\ph}{{}\,\widehat{\!p}{}}
\newcommand{\muh}{{}\,\widehat{\!\mu}{}}
\newcommand{\jo}{{}\,\overline{\!j}{}}
\newcommand{\Lo}{{}\,\overline{\!L}{}}
\newcommand{\XXo}{{}\,\overline{\!X}{}}
\newcommand{\Xo}{{}\,\overline{\!\mathcal X}{}}
\newcommand{\Yo}{{}\,\overline{\!\mathcal Y}{}}

\newcommand{\dd}{\partial}
\newcommand{\ddd}{{\rm d}}
\newcommand{\ee}{{\bf e}}
\newcommand{\Ff}{F_{\!f}}
\newcommand{\Gp}{\Gamma_{\!+}}
\newcommand{\Gf}{\Gamma_{\!f}}
\newcommand{\Gg}{\Gamma_{\!g}}
\newcommand{\Spx}{{\rm Supp}_{(x)}}
\newcommand{\Sp}{{\rm Sp}}
\newcommand{\SpW}{{}^W\!\1{\rm Sp}}
\newcommand{\CFf}{{\rm CF}_{\hskip-2pt f}}
\newcommand{\CFfb}{{\rm CF}_{\hskip-2pt f,{\bf b}}}
\newcommand{\CFfsb}{{\rm CF}_{\hskip-2pt f,{\bf sb}}}
\newcommand{\CFfdb}{{\rm CF}_{\hskip-2pt f,{\bf db}}}
\newcommand{\CFfnb}{{\rm CF}_{\hskip-2pt f,{\bf nb}}}
\newcommand{\CFfin}{{\rm CF}_{\hskip-2pt f,{\rm in}}}
\newcommand{\CFfla}{{\rm CF}_{\hskip-2pt f,\la}}
\newcommand{\CFfinla}{{\rm CF}_{\hskip-2pt f,{\rm in},\la}}

\newcommand{\eq}{\,{=}\,}

\newcommand{\pl}{\1{+}\1}
\newcommand{\mi}{\1{-}\1}
\newcommand{\bl}{\bigl}
\newcommand{\br}{\bigr}

\newcommand{\ssb}{\raise.15ex\h{${\scriptscriptstyle\bullet}$}}
\newcommand{\ssc}{\,\raise.15ex\h{${\scriptstyle\circ}$}\,}
\newcommand{\onto}{\twoheadrightarrow}
\newcommand{\into}{\hookrightarrow}
\newcommand{\simto}{\,\,\rlap{\hskip1.5mm\raise1.4mm\hbox{$\sim$}}\hbox{$\longrightarrow$}\,\,}
\newcommand{\nsset}{\rlap{\raise.5mm\h{$\subset$}}\raise-1.5mm\h{$\rlap{\raise.5mm\h{$\scriptscriptstyle\,\,/$}}-$}}
\begin{document}
\title[Descent of nearby cycle formula]
{Descent of nearby cycle formula\\for Newton non-degenerate functions}
\author[M. Saito]{Morihiko Saito}
\address{M. Saito : RIMS Kyoto University, Kyoto 606-8502 Japan}
\email{msaito@kurims.kyoto-u.ac.jp}
\begin{abstract} We prove a descent theorem of nearby cycle formula for Newton non-degenerate functions at the origin as well as its motivic version (without assuming the convenience condition). This is used in some papers without any proof although its proof is quite nontrivial because of the existence of coordinate hyperplanes, which is completely neglected in the literature about the descent theorem. In the isolated singularity case, it implies some well-known formula for the number of Jordan blocks of the Milnor monodromy with the theoretically maximal size, using a standard estimate of weights. It also provides a proof of a modified version of the Steenbrink conjecture on spectral pairs for non-degenerate functions with simplicial Newton polytopes in the isolated singularity case (which is false in the non-simplicial case).
\end{abstract}
\maketitle
\centerline{\bf Introduction}
\bsn
Let $f$ be a Newton non-degenerate holomorphic function of $n$ variables with $\Gp(f)$ the Newton polytope at the origin, fixing local coordinates $x_1,\dots,x_n$. We have the dual fan $\Si$ in $(\R_{\ges 0})^n$ (see \cite[9.1]{Va1}) and also a toric variety $\X$ with a proper morphism
$$\X\buildrel{\pi}\over\longrightarrow X:=\C^n,
$$
\par\nin inducing an isomorphism over $X^*:=(\C^*)^n$. Here we use the {\it associated analytic spaces\1} rather than algebraic varieties (since we have to deal with {\it local systems}).
\sk
Let $\Gf$ be the union of {\it compact\1} faces of $\Gp(f)$. There is a stratification
$$\X_0^{\rm ex}:=\pi^{-1}(0)=\msqcup_{\si\les\Gf}\,\X_{\si},$$
\par\nin where $\X_{\si}\cong(\C^*)^{d_{\si}}$ for a face $\si\les\Gf$ with $d_{\si}:=\dim\si$. These $\X_{\si}$ are {\it orbits\1} of the natural torus action on $\X$. This stratification is {\it finer\1} than the one used for the calculation of motivic nearby fibers in \cite[3.3]{DL}, since it gives a stratification of the inverse image of the {\it union\1} of $f^{-1}(0)$ and the {\it coordinate hyperplanes.} The codimension of each stratum must be {\it adjusted\1} by subtracting the number of the proper transforms of coordinate hyperplanes containing the stratum, and this is closely related to the numbers $k(\si)$ explained below. These do not seem to be mentioned in the literature.
\sk
Let $\X_0^{\rm pr}\subset\X$ be the proper transform of $f^{-1}(0)$. This is a hypersurface of some {\it open neighborhood\1} of $\X_0^{\rm ex}\subset\X$ in the classical topology.
Set $\X_{\si}^{\circ}:=\X_{\si}\setminus\X_0^{\rm pr}$ with $j_{\si}:\X_{\si}^{\circ}\into\X_{\si}$ the inclusion.
\sk
Let $K_0^{T_s}({\rm MHS})$ be the Grothendieck ring of mixed $\Q$-Hodge structures endowed with an endomorphism $T_s$ of finite order. Denoting the Milnor fiber of $f$ around 0 by $\Ff$, we have the Euler characteristic
$$[\1\chi(\Ff,\Q),T_s]\,\in\,K_0^{T_s}({\rm MHS}),$$
\par\nin where $T_s$ is the semisimple part of the {\it Jordan decomposition\1} $T=T_sT_u$ of the monodromy $T$. (Note that the action of the monodromy $T$ on the $H^k(\Ff,\Q)$ is {\it not\1} an endomorphism of mixed Hodge structure, unless $T$ is semisimple, that is, $T=T_s$.) Set
$$\theta:=[\1\Q(-1),{\rm id}\1]\,\in\,K_0^{T_s}({\rm MHS}),$$
\par\nin which is the class of a mixed $\Q$-Hodge structure of type $(1,1)$ with trivial action of $T_s$.
\sk
Let $x_1,\dots,x_n$ be the coordinates used for the Newton polytope. Put
$$J_f:=\bl\{i\in[1,n]\,\big|\,\{x_i\eq 0\}\subset\{f\eq 0\}\br\}.$$
\par\nin We have $J_f=\emptyset$ if $f$ is convenient (or $f$ has an isolated singularity with $n\ges 3$).
For a face $\si\les\Gf$, set $d(\si):=\dim C(\si)\,(=d_{\si}{+}1)$ with $C(\si)\subset\R^n$ the cone of $\si$, and
$$k(\si):=\min\bl\{\1|I|\,\1\big|\1\,J_f\subset I\subset\{1,\dots,n\},\,\,\R^{I}\supset\si\br\},$$
\par\nin where $\R^I:=\mcap_{i\notin I}\,\{x_i=0\}\subset\R^n$.
Let $\de_{\si}$ be the positive integer satisfying
\htt{1}{}
$$\ell_{\si}(V(\si)\cap\Z^n)=\tfrac{\!\!1}{\de_{\si}}\1\Z,
\leqno(1)$$
\par\nin where $V(\si)\subset\R^n$ is the vector subspace spanned by $\si$, and $\ell_{\si}$ is the linear function on $V(\si)$ with $\ell_{\si}^{-1}(1)\supset\si$. In this paper we show the following.
\par\htt{T1}{}\msn
{\bf Theorem\,\,1.} {\it For a Newton non-degenerate function $f$, there is an equality in $K_0^{T_s}({\rm MHS})$
\htt{2}{}
$$\aligned\bl[\bl(\chi(\Ff,\Q),T_s\br)\br]={}&\msum_{\si\les\Gf}\,(-1)^{d(\si)}(1-\theta)^{k(\si)-d(\si)}\1\psi_{\si}\,,\\
\h{with}\,\,\,\q\q\q\psi_{\si}:={}&\bl[H^{d_{\si}}_c\bl(\X_{\si},\RR(j_{\si})_*(L_{\si},T_s)\br)\br],\endaligned
\leqno(2)$$
\par\nin Here $\si$ runs over all the faces of $\Gf$ $($not including $\,\emptyset)$, and $L_{\si}$ is a $\Q$-local system of rank $\de_{\si}$ on $\X_{\si}^{\circ}$, which is identified with a variation of Hodge structure of type $(0,0)$, and is endowed with an action of $T_s$ of finite order. Moreover ${\rm Ker}(T_s{-}\la)$ in the complex scalar extension $\C\otimes_{\Q}L_{\si}$ for $\la\in\C^*$ is a local system of rank~$1$ if $\la^{\de_{\si}}=1$, and vanishes otherwise.}
\ms
This may be called a {\it descent theorem\1} of nearby cycle formula for Newton non-degenerate functions, since a similar formula is known for a desingularization by a smooth subdivision of $\Si$ (where the cohomology with compact supports in the definition of $\psi_{\si}$ is replaced by the Euler characteristic with compact supports up to sign depending of $d_{\si}$), see Theorem\,\,\hl{TA.2}{A.2} below for a {\it motivic\1} version. The local systems $L_{\si}$ can be defined in a similar way to motivic nearby fibers of Denef and Loeser \cite[3.3]{DL}.
\sk
Here a problem is that the $k(\si)$ (or more precisely $\kc_{\eta_{\si}}$ in (\hl{2.3.2}{2.3.2}) below) are {\it unstable\1} under the morphism of toric orbits, and this makes the proof of Theorem\,\,\hl{T1}{1} quite complicated, see the end of the proof of Theorem\,\,\hl{T1}{1} in \hl{2.3}{2.3} below. This comes from the existence of coordinate hyperplanes explained above. There is also a (minor) problem of an ``obvious combinatorial identity" used without any proof, which is not easy to prove {\it combinatorially\1} without applying a ``topological method" (employing Euler numbers) as in the proof of Proposition\,\,\hl{P2.3}{2.3} below. These points do not seem to be studied carefully in the papers quoted by \cite{Sta} as relevant results. Except for the argument related to the problem of $k(\si)$, the proof of Theorem\,\,\hl{T1}{1} is rather an easy consequence of Remarks\,\,\hl{R2.2c}{2.2c}--\hl{R2.2e}{e} below. (It is also possible to prove Theorem\,\,\hl{T1}{1} by induction on $\dim\eta_{\si}-\dim\xi$ as in the proof of Theorem\,\,\hl{TA.2}{A.2} below.)
\sk
We say that $f$ is {\it simplicial,} if every compact face of $\Gp(f)$ is a simplex. In this case, we can define a graded $\C$-algebra $B_{\si}$ to be the $\C$-vector space spanned by the monomials $x^{\nu}$ with $\nu\in\N^n$ contained in the interior $E_{\si}^{\circ}$ of the {\it parallelotope\1} $E_{\si}$ spanned by the vertices $v_j$ of a face $\si\les\Gf$, that is, $E_{\si}^{\circ}=\bl\{\msum_j\,r_jv_j\mid r_j\in(0,1)\br\}$. Here the grading is defined so that any vertex $v_j$ of $\si$ has degree 1. Let $q_{\si}(t)$ be the Poincar\'e polynomial of $B_{\si}$, see also (\hl{1.4.1}{1.4.1}--\hl{1.4.2}{2}) below. For $\tau\les\Gf$, set
$$\qh_{\tau}(t)=\msum_{\emptyset\les\si\les\tau}\,q_{\si}(t).$$
\par\nin Here $\si$ runs over any face of $\tau$ {\it including\1} $\emptyset$, and $q_{\emptyset}(t):=1$.
\sk
Let $\Sp'_f(t)$ be the {\it Hodge spectrum\1} of $f$, see \cite{ste}, \cite{DL}, and also \hl{1.2}{1.2} below. This is essentially the {\it dual\1} of the usual Steenbrink spectrum \cite{St2} (up to the shift of $t$ by 1) in the {\it non-isolated\1} singularity case, see also \cite[1.8]{JKSY1}, etc. Theorem\,\,\hl{T1}{1} has the following.
\par\htt{C1}{}\msn\vbox{\nin
{\bf Corollary\,\,1.} {\it If $f$ is Newton non-degenerate and simplicial, then we have the equality
\htt{3}{}
$$\Sp'_f(t)=\msum_{\emptyset\les\si\les\Gf}\,(-1)^{n-d(\si)}\,(1{-}t)^{k(\si)-d(\si)}\1\qh_{\si}(t),
\leqno(3)$$
\par\nin with $\si$ running over any faces of $\Gf$ {\it including\1} $\emptyset$, where $k(\emptyset)=d(\emptyset)=0$.}}
\ms
Note that reduced cohomology is used for the definition of spectrum, and the summand for $\si=\emptyset$ corresponds to the difference between the usual and reduced cohomologies.
We can show that Corollary\,\,\hl{C1}{1} implies \cite[Theorem 2]{JKSY2}, although this requires some non-trivial assertion on combinatorics of triangulations of triangles, see Remarks\,\,\hl{R2.4e}{2.4e}--\hl{R2.4f}{f} below.
\sk
In the {\it convenient\1} Newton non-degenerate case (where the singularities are isolated \cite{Ko}), the equality (\hl{3}{3}) follows from \cite[Theorem 5.7]{St1} (see also \cite[1.5]{JKSY2}) using \cite{exp}. We say that $f$ is {\it convenient,} if $\Gp(f)$ intersects any coordinate axis of $\R^n$. Under this hypothesis, it is {\it unnecessary\1} to assume that $f$ is simplicial, see (\hl{4}{4}) and also {1.4}{1.4} below.
\sk
Using \cite{mic}, \cite{bfun}, \cite{hi}, we can deduce from Corollary\,\,\hl{C1}{1} the following.
\par\htt{C2}{}\msn
{\bf Corollary\,\,2.} {\it Assume $f$ is Newton non-degenerate and simplicial. Let $c\in\Q_{>0}$ such that $(c,\dots,c)\in\dd\Gp(f)$. If $c>1$, then $-1/c$ is the maximal root of the Bernstein-Sato polynomial $b_{f,0}(s)$.}
\ms
In the isolated singularity case, this is a consequence of \cite{EL}, \cite{Ma}, \cite{exp}, \cite{SS}, \cite{Va2}, where the minimal spectral number coincides with the maximal root of $b_{f,0}(s)/(s{+}1)$ up to sign; hence Corollary\,\,\hl{C2}{2} holds also in the case $c\les 1$ by replacing $b_{f,0}(s)$ with $b_{f,0}(s)/(s{+}1)$. (It has been informed that some argument in \cite{EL} does not work when $c<1$.)
\sk
From now on, we assume $f$ has an {\it isolated singularity\1} at 0 in this introduction. (However, we do not assume $f$ is convenient.)
Let $\CFf^k$ be the set of $k$-dimensional compact faces of $\Gp(f)$, and $\CFfin^k$ be the subset consisting of internal faces. (A face is called {\it internal\1} if it is not contained in any coordinate hyperplane of $\R^n$.) For $\la\in\C^*$, $k\in\N$, set
$$\aligned\CFfla^k&:= \{\si\in\CFf^k\mid\la^{\de_{\si}}=1\},\\ \CFfinla^k&:=\CFfin^k\cap\CFfla^k,\endaligned$$
\par\nin \vskip-2mm\nin
and
\vskip-4mm
$$\begin{array}{ll} \,l_{\si}:=|\Z^n\cap\si|-1&\bl(\si\in\CFf^1\br),\\ \be_{\si}:=\#\bl\{\tau\in\CFfin^1\,\big|\,\tau\ges\si\br\}&\bl(\si\in\CFf^0).\raise5mm\h{}\end{array}$$
\par\nin Let $n_{\la,k}$ be the number of Jordan blocks of the monodromy $T$ for eigenvalue $\la$ with size $k$.
As a consequence of Theorem\,\,\hl{T1}{1}, we get a proof of the following formula which seems to be known to certain specialists (see for instance \cite{Sta} and references there).
\par\htt{C3}{}\msn
{\bf Corollary\,\,3.} {\it Assume $f$ is Newton non-degenerate, and has an isolated singularity at the origin. Then we have the equalities}
$$n_{\la,k}=\begin{cases}\,\bl|\1\CFfinla^0\br|&(\la\,{\ne}\,1,\,k\eq n),\\ \,\msum_{\si\in\CFfinla^1}\,l_{\si}\,-\,\msum_{\si\in\CFfla^0}\,\be_{\si}\raise5mm\h{}&(\la\,{\ne}\,1,\,k\eq n{-}1),\\ \,\bl|\1\Z_{>0}^n\cap\1\mcup_{\si\in\CFfin^1}\1\si\1\br|\raise5mm\h{}&(\la\eq 1,\,k\eq n{-}1).\end{cases}$$
\par\nin \ms
Note that $n_{1,n}=0$ as a consequence of (\hl{1.2.4}{1.2.4}) below.
This corollary immediately follows from Theorem\,\,\hl{T1}{1} using a {\it standard estimate of weights\1} as in \cite[4.5.2]{mhm} (especially for $n_{\la,n}$). Here we do not need a theory of the duals of logarithmic complexes on toric varieties. (The latter theory could be justified partially by using \cite[3.11]{mhm} for instance.)
The above formula for $\la\,{\ne}\,1$, $k\eq n{-}1$ is equivalent to a well-known one (see for instance \cite{Sta} and also Remark\,\,\hl{R2.5b}{2.5b} below).
\sk
Assume now $f$ is furthermore {\it simplicial}. In this case there is no counter-example to a modified version of Steenbrink's conjecture on spectral pairs of Newton non-degenerate functions, see \cite{Da2} (and Remark\,\,\hl{R1.4a}{1.4a} below). Let $(\al_i,w_i)$ ($i\in[1,\,\mu_f])$ be the {\it spectral pairs\1} of $f$, where $\mu_f$ is the Milnor number, the $\al_i\in\Q$ are {\it spectral numbers\1}, and the $w_i\in\N$ are {\it modified weights.} The definition in this paper essentially coincides with the one in \cite{St1}, and is different from the one in \cite{SSS} by the shift of spectral numbers $\al_i$ by $1$ and the change of modified weights $w_i$ by the involution of $\Z$ defined by $w\mapsto 2n\mi2\mi w$, see \hl{1.2}{1.2} below.
\sk
By \cite[Theorem 5.7]{St1} (see also \cite[1.5]{JKSY2}) and using \cite{exp}, we have the equalities
\htt{4}{}
$$\aligned\Sp_f(t):=\msum_{i=1}^{\mu_f}\,t^{\al_i}={}&\msum_{\emptyset\les\si\les\Gf}\,(-1)^{n-d(\si)}\,(1{-}t)^{k(\si)-d(\si)}\1\qh_{\si}(t)\\ ={}&\msum_{\emptyset\les\si\les\Gf}\,r_{\si}(t)q_{\si}(t)\q\q\q\q\h{with}\\ r_{\si}(t):={}&\msum_{\si\les\tau\les\Gf}\,(-1)^{n-k(\tau)}\,(t{-}1)^{k(\tau)-d(\tau)}.\endaligned
\leqno(4)$$
\par\nin The summations $\msum_{\emptyset\les\si\les\Gf}$ in (\hl{4}{4}) are taken over any faces of $\Gf$ {\it including\1} $\emptyset$. We call $r_{\si}(t)$ the {\it combinatorial polynomial}. The assertion (\hl{4}{4}) holds {\it without\1} assuming $f$ is simplicial. It was originally proved in the $f$ {\it convenient\1} case. However, it holds also in the non-convenient simplicial case by (\hl{3}{3}) assuming $f$ has an isolated singularity. Observe the striking similarity between the right-hand sides of (\hl{2}{2}) and (\hl{4}{4}), although these come from totally different methods.
\sk
Using Theorem\,\,\hl{T1}{1}, we can prove a modified Steenbrink conjecture (see \cite{St1}) as follows.
\par\htt{T2}{}\msn
{\bf Theorem\,\,2.} {\it Assume $f$ is Newton non-degenerate, simplicial, and has an isolated singularity at the origin. Then the generating function of the spectral pairs, which is called the weighted spectrum, is given by
\htt{5}{}
$$\SpW_f(t,u):=\msum_{i=1}^{\mu_f}\,t^{\al_i}u^{w_i}=\msum_{\emptyset\les\si\les\Gf}\,r_{\si}(tu^2)\1 u^{d_{\si}}\1q_{\si}(t),
\leqno(5)$$
\par\nin where $\si$ runs over any faces of $\Gf$ including $\emptyset$ with $d_{\emptyset}:=-1$.}
\ms
We have the {\it symmetry\1} of combinatorial polynomials $r_{\si}(t)$ in the convenient simplicial case, see \cite[Theorem A]{JKSY2}. This is compatible with the symmetry of spectral pairs for the weight filtration (\hl{1.2.5}{1.2.5}). In the proof of Theorem\,\,\hl{T2}{2}, it is better to adopt the definition of spectrum like {\it Hodge spectrum\1} in Corollary\,\,\hl{C1}{1} (using the symmetry of spectral numbers in the isolated singularity case), see \hl{1.2}{1.2} below.
In the {\it non-simplicial\1} case, it does not seem easy to determine the equivariant Euler characteristic Hodge numbers of the vanishing cohomology of $f$ (where equivariant means that the Hodge numbers are given for each monodromy eigenvalue). It is unclear whether there is really an algorithm to determine them except for the case where all the compact faces of the Newton polytope are {\it prime\1}, that is, its dual fan is simplicial. Some arguments in papers quoted by \cite{Sta} as relevant results may work only in the prime case, see also Remark\,\,\hl{RA.2b}{A.2b} below.
\sk
This work was partially supported by JSPS Kakenhi 15K04816.
\sk
In Section 1 we recall basics of spectrum and spectral pairs.
In Section 2 we review some toric geometry, and prove the main theorems.
In Appendix we show the descent theorem for motivic nearby fibers.
\msn
{\bf Convention.} For a subset $\ga\subset\R^n$, we denote respectively by $A(\ga)$, $V(\ga)$ the smallest affine and vector subspaces containing $\ga$ in this paper. (This is different from \cite{JKSY1}, \cite{JKSY2}, where $A_{\ga}$, $V_{\ga}$ are used respectively. In this paper $A_{\ga}$ means a $\C$-algebra when $\ga=\si$.)
\sk
We denote by $\de_{\ga}$ the largest {\it non-negative\1} integer $d$ such that $d\1\nu\in A(\ga)$ for some $\nu\in\Z^n$. This is compatible with a definition after (\hl{1}{1}). We can also define $\de_{\ga}$ using a linear function $\ell_{\ga}$ on $V(\ga)$ as in (\hl{1}{1}) when $0\notin A(\ga)$. Note that $\de_{\ga}=0$ if and only if $0\in A(\ga)$.
\bs\ms\vbox{
\centerline{Contents}
\ms
\par\hl{S1}{1. Spectrum and spectral pairs}\hfill 5
\par\q\hl{1.1}{1.1.~Spectrum}\hfill 5
\par\q\hl{1.2}{1.2.~Spectral pairs}\hfill 5
\par\q\hl{1.3}{1.3.~Newton non-degenerate case}\hfill 6
\par\q\hl{1.4}{1.4.~Steenbrink formula for spectrum}\hfill 7
\par\hl{S2}{2. Relation to the toric geometry}\hfill 8
\par\q\hl{2.1}{2.1.~Toric varieties}\hfill 8
\par\q\hl{2.2}{2.2.~Construction of nearby cycle sheaves}\hfill 12
\par\q\hl{2.3}{2.3.~Proof of Theorem}\,\,\hl{T1}{1}\hfill 14
\par\q\hl{2.4}{2.4.~Proofs of Corollary}\,\,\hl{C1}{1} \hl{2.4}{and Theorem}\,\,\hl{T2}{2}\hfill 18
\par\q\hl{2.5}{2.5.~Proof of Corollary}\,\,\hl{C3}{3}\hfill 22
\par\q\hl{2.6}{2.6.~Proof of Corollary}\,\,\hl{C2}{2}\hfill 22
\par\hl{A}{Appendix. Descent theorem for motivic nearby fibers}\hfill 23
\par\q\hl{A.1}{A.1.~Relative Grothendieck ring}\hfill 23
\par\q\hl{A.2}{A.2.~Motivic descent theorem}\hfill 24}
\bsn\msn\htt{S1}{}\par
\vbox{\centerline{\bf 1. Spectrum and spectral pairs}
\bsn
In this section, we recall basics of spectrum and spectral pairs.}
\par\htt{1.1}{}\msn
{\bf 1.1.~Spectrum.} Let $f:(\C^n.0)\to(\C,0)$ be a holomorphic function. We denote the Milnor fiber of $f$ by $\Ff$. We have the canonical mixed Hodge structure on the vanishing cohomology $\Ht^j(\Ff,\Q)$ (using for instance mixed Hodge modules \cite{mhm}), where $\Ht$ denotes the reduced cohomology. Let $F,W$ be the Hodge and weight filtrations.
\sk
Let $T=T_sT_u$ be the {\it Jordan decomposition\1} of the monodromy $T$. Note that this $T$ is the {\it inverse\1} of the Milnor monodromy. This is closely related to some confusion in the definition of spectrum in \cite{St1}, see for instance \cite{DS}. Set
$$\Ht^j(\Ff,\C)_{\la}:={\rm Ker}(T_s\mi\la)\,\subset\,\Ht^j(\Ff,\C).$$
\par\nin We define the {\it Hodge spectrum\1} $\Sp'_f(t)=\msum_{\al\in\Q}\,m'_{f,\al}\1t\1^{\al}$ by
\htt{1.1.1}{}
$$\aligned&m'_{f,\al}=\msum_{j=0}^{n-1}\,(-1)^{n-1-j}\dim_{\C}\Gr_F^p\Ht^j(\Ff,\C)_{\la}\\ &\h{with}\q\q\q p=[\al],\q\la=\exp(2\pi\sqrt{-1}\al),\endaligned
\leqno(1.1.1)$$
\par\nin see \cite{ste}, \cite{DL}, etc. This is essentially the {\it dual\1} of the Steenbrink spectrum $\Sp_f(t)$, that is,
\htt{1.1.2}{}
$$\Sp'_f(t)=t^n\1\Sp_f(1/t).
\leqno(1.1.2)$$
\par\nin Indeed, the conditions $p\eq[\al]$, $\la\eq e^{2\pi\sqrt{-1}\al}$ are replaced respectively with $p\eq[n{-}\al]$ and $\la\eq e^{-2\pi\sqrt{-1}\al}$ in the usual definition of spectrum, see for instance \cite{BS1}, \cite{BS2}, \cite{DS}, \cite{JKSY1}, etc.
Note that the Steenbrink spectrum in \cite{St2} is shifted by $-1$ so that the spectral numbers (that is, the $\al$ with $m'_{f,\al}\ne 0$) are contained in $(-1,n{-}1)$ instead of $(0,n)$.
\par\htt{1.2}{}\msn
{\bf 1.2.~Spectral pairs.} In the isolated singularity case ($n\ges 2$), we define the {\it spectral pairs\1} $(\al_i,w_i)$ ($i\in[1,\,\mu_f])$ in this paper as follows (see also \cite{St1}, \cite[2.1--2]{ste}, \cite{JKSY2}):
\htt{1.2.1}{}
$$\aligned&\#\{\1i\mid(\al_i,w_i)\eq(\al,w)\}=\dim\Gr_F^p\Gr^W_{w+\de_{\la,1}}H^{n-1}(\Ff,\C)_{\la}\\ &\q\q\h{with}\q\q\q p=[\al],\q\la=\exp(2\pi\sqrt{-1}\al),\endaligned
\leqno(1.2.1)$$
\par\nin where $\mu_f$ is the Milnor number (and $\de_{\la,1}=1$ if $\la=1$, and 0 otherwise). This definition of the $\al_i$ is different from the usual one as is explained after (\hl{1.1.1}{1.1.1}). However, it does not cause a problem because of the {\it symmetry\1} explained just below.
The numbers $\al_i\in\Q$ and $w_i\in\N$ are respectively called the {\it spectral numbers\1} and the {\it modified weights.}
The generating polynomial of the spectral pairs is called the {\it weighted spectral,} and is denoted by $\SpW_f(t,u)$, that is,
$$\SpW_f(t,u)=\msum_{i=1}^{\mu_f}\,t^{\al_i}u^{w_i}.$$
\par\nin There are equalities
\htt{1.2.2}{}
$$\Sp_f(t)=\msum_{i=1}^{\mu_f}\,t^{\al_i}=\SpW(t,1).
\leqno(1.2.2)$$
\par\nin Here we use the {\it self-duality\1} (or {\it symmetry})
\htt{1.2.3}{}
$$\Sp_f(t)=\Sp_f(t^{-1})\1t^n,\q\h{that is,}\q\al_i\eq\al_j\,\,\,\,(\1i{+}j\eq\mu_f{+}1),
\leqno(1.2.3)$$
\par\nin which follows from the {\it Hodge symmetry\1} of the graded pieces of the vanishing cohomology $\Gr^W_kH^{n-1}(\Ff,\C)$ together with the {\it monodromical property\1} of the weight filtration $W$ as in (\hl{1.2.4}{1.2.4}) below (using also the assertion that $T_s$ is {\it defined over\1} $\R$).
\par\htt{R1.2a}{}\msn
{\bf Remark\,\,1.2a.} The above definition of spectral pairs coincides essentially with the one in \cite{St1}, and is different from the one in \cite{SSS} as is explained before (\hl{4}{4}) in the introduction. The generating polynomial for the spectral pairs in the sense of \h{\it loc.\,cit.} is given by
$$\SpW_f(t,u^{-1})\1t^{-1}\1 u^{2n-2}.$$
\par\nin In the case $f=x^p\pl y^p\pl z^p\pl xyz\q(p>3)$, for instance, the generating polynomials are respectively given by
$$\aligned tu\pl t^2u^3\pl&3\,\msum_{k=1}^{p-1}\,t^{\,k/p+1}u^2,\\ u^3\pl tu\pl&3\,\msum_{k=1}^{\,p-1}\,t^{k/p}u^2.\endaligned$$
\par\nin \par\htt{R1.2b}{}\msn
{\bf Remark\,\,1.2b.} We have the decomposition by {\it unipotent\1} and {\it non-unipotent\1} monodromy part of the vanishing cohomology:
$$H^{n-1}(\Ff,\C)=H^{n-1}(\Ff,\C)_1\oplus H^{n-1}(\Ff,\C)_{\ne 1},$$
\par\nin where the last term is defined by
$$H^{n-1}(\Ff,\C)_{\ne 1}:=\mopl_{\la\ne 1}\,H^{n-1}(\Ff,\C)_{\la}.$$
\par\nin The weight filtration $W$ on $H^{n-1}(\Ff,\C)_{\ne 1}$ and $H^{n-1}(\Ff,\C)_1$ {\it coincides\1} with the {\it monodromy filtration\1} shifted respectively by $n{-}1$ and $n$, that is,
\htt{1.2.4}{}
$$\aligned N^j:\Gr_{n-1+j}^WH^{n-1}(\Ff,\C)_{\ne 1}&\simto\Gr_{n-1-j}^WH^{n-1}(\Ff,\C)_{\ne 1}\q(j\in\N),\\ N^j:\Gr_{n+j}^WH^{n-1}(\Ff,\C)_1&\simto\Gr_{n-j}^WH^{n-1}(\Ff,\C)_1\q(j\in\N),\endaligned
\leqno(1.2.4)$$
\par\nin where $N=\log T_u$. The proof of (\hl{1.2.4}{1.2.4}) is highly non-trivial, and some argument as in \cite[4.2.2]{mhp} must be needed, since the weight filtration $W$ is {\it not\1} defined as a shifted monodromy filtration in \cite{St1} and the passage from the $E_1$-term to the $E_2$-term is quite non-trivial.
\par\htt{R1.2c}{}\msn
{\bf Remark\,\,1.2c.} The spectral pairs $(\al_i,w_i)$ of $f$ are equivalent to the {\it weight decomposition}
$$\Sp_f(t)=\msum_{w=0}^{2n-2}\,\Sp_f(t)_{(w)}\q\h{with}\q\Sp_f(t)_{(w)}=\msum_{w_i=w}\,t^{\al_i}.$$
\par\nin From (\hl{1.2.4}{1.2.4}) we can deduce the {\it symmetry\1} of the weight decomposition
\htt{1.2.5}{}
$$t^j\1\Sp_f(t)_{(n-1-j)}=\Sp_f(t)_{(n-1+j)}\q(j\in\N).
\leqno(1.2.5)$$
\par\nin \par\htt{1.3}{}\msn
{\bf 1.3.~Newton non-degenerate case.} We denote by $\Gp(f)$ the {\it Newton polytope\1} of $f$. This is the convex hull of the union of $\nu+\R_{\ges 0}^n$ for $\nu\in\Spx\1f$ with
\htt{1.3.1}{}
$$\Spx\1f:=\{\nu\in\N^n\mid a_{\nu}\ne 0\1\}\q\h{for}\q f=\msum_{\nu}\,a_{\nu}x^{\nu}\in\C\{x\},
\leqno(1.3.1)$$
\par\nin where $x_1,\dots,x_n$ are the coordinates of $\C^n$.
\sk
We say that $f$ is ({\it Newton\1}) {\it non-degenerate,} or more precisely, $f$ has {\it non-degenerate Newton boundary,} if we have for any {\it compact\1} face $\si\subset\Gp(f)$
\htt{1.3.2}{}
$$\mcap_{i=1}^n\,\bl\{\1 x_i\dd_{x_i}f_{\si}=0\br\}\cap(\C^*)^n=\emptyset,
\leqno(1.3.2)$$
\par\nin where $f_{\si}:=\mopl_{\nu\in\si}\,a_{\nu}x^{\nu}$ with $a_{\nu}$ as in (\hl{1.3.1}{1.3.1}), see \cite{Ko}, \cite{Va1}, etc.
\sk
Assume now $f$ is non-degenerate and moreover {\it convenient\1} (that is, $\Gp(f)$ intersects every coordinate axis of $\R^n$). These conditions imply that $f$ has an isolated singularity at 0, see \cite{Ko}. For $h\in\C\{x\}$, set
\htt{1.3.3}{}
$$v_f(h):=\max\bl\{r\in\R\mid\ob\pl\Spx\1h\,\subset\,r\,\Gp(f)\br\},
\leqno(1.3.3)$$
\par\nin with $\ob:=(1,\dots,1)$. The minimal spectral number $\al_1$ coincides with $v_f(1)$. More generally, the $V$-filtration on $\C\{x\}/(\dd f)$ induced from the Brieskorn lattice (see for instance \cite{SS}) coincides with the one induced from the Newton filtration $V_N^{\ssb}$ on $\C\{x\}$ defined by
\htt{1.3.4}{}
$$V_N^{\al}\C\{x\}:=\bl\{h\in\C\{x\}\mid v_f(h)\ges\al\br\}\q(\al\in\Q),
\leqno(1.3.4)$$
\par\nin see \cite{exp}.We get for $\al\in\Q$
\htt{1.3.5}{}
$$\#\bl\{i\in[1,\mu_f]\,\big|\,\al_i=\al\br\}=\dim_{\C}\Gr_{V_N}^{\al}\bl(\C\{x\}/(\dd f)\br),
\leqno(1.3.5)$$
\par\nin where $(\dd f)\subset\C\{x\}$ denotes the Jacobian ideal.
\sk
We define $v'_f$, $V'_N$ in the same way as in (\hl{1.3.3}{1.3.3}--\hl{1.3.4}{4}) except that $\ob+$ is omitted in (\hl{1.3.3}{1.3.3}). We have $v_f(h)=v'_f(x^{\ob}h)$ with $x^{\ob}=x_1\cdots x_n$.
\par\htt{1.4}{}\msn
{\bf 1.4.~Steenbrink formula for spectrum.} Let $f$ be a holomorphic function with convenient non-degenerate Newton boundary (in particular, $f$ has an isolated singularity at 0), see \hl{1.3}{1.3}. We do {\it not\1} assume $f$ is simplicial as in the introduction. We have the graded $\C$-algebra
$$A_{\Ga}:=\Gr_{V'_N}^{\ssb}\C\{x\}.$$
\par\nin Let $\Gf$ be the union of {\it compact\1} faces of $\Gp(f)$. For $\si\les\Gf$, let $A_{\si}$ be the graded subalgebra of $A_{\Ga}$ generated by the classes of monomials $[x^{\nu}]$ with $\nu$ contained in the cone $C(\si)$ of $\si$ in $\R^n$. Let $p_{\si}(t)$ be the Hilbert-Poincar\'e series of the graded vector space $A_{\si}$ (which is a fractional power series). Set
\htt{1.4.1}{}
$$\aligned\qh_{\si}(t)&:=(1-t)^{d(\si)}p_{\si}(t)\q\h{with}\q d(\si):=\dim C(\si),\\ q_{\si}(t)&:=\msum_{\emptyset\les\tau\les\si}\,(-1)^{d(\si)-d(\tau)}\qh_{\tau}(t).\endaligned
\leqno(1.4.1)$$
\par\nin with summation taken over any faces $\tau$ of $\si$ including $\si$ and $\emptyset$. Here $C(\emptyset)=\{0\}$ and $\qh_{\emptyset}(t)=q_{\emptyset}(t)=1$. We have $d(\si)=d_{\si}{+}1$ (even if $\si=\emptyset$). Note that $\qh_{\si}(t)$ is a fractional power polynomial, see \cite{Ko}, \cite{St2}. It is known ({\it loc.\,cit.}) that
\htt{1.4.2}{}
$$\qh_{\si}(t)=\msum_{\emptyset\les\tau\les\si}\,q_{\tau}(t).
\leqno(1.4.2)$$
\par\nin (Note that (\hl{1.4.2}{1.4.2}) is trivial in the simplicial case.) This definition of $q_{\si}$ is compatible with the one in the introduction when $f$ is {\it simplicial}. Set as in the introduction
$$k(\si)=\min\{|I|\mid C(\si)\subset\R^I\}\q\h{with}\q\R^I:=\mcap_{i\notin I}\,\{x_i=0\}.$$
\par\nin \sk
Consider the graded $A_{\Ga}$-modules
$$\aligned B_{\Ga}&:=\Gr_{V_N}^{\ssb}\bl(\C\{x\}/(\dd f)\br),\\ B'_{\Ga}&:=\Gr_{V'_N}^{\ssb}\bl(\C\{x\}/(x_1f_1,\dots,x_nf_n)\br),\endaligned$$
\par\nin where $f_i:=\dd_{x_i}f$. The Poincar\'e polynomial of $B_{\Ga}$ can be described as
\htt{1.4.3}{}
$$p_{B_{\Ga}}(t)=\msum_{\emptyset\les\si\les\Gf}\,(-1)^{n-d(\si)}\,(1{-}t)^{k(\si)-d(\si)}\1\qh_{\si}(t),
\leqno(1.4.3)$$
\par\nin where the summation is taken over any faces of $\Gf$ including $\emptyset$, see \cite[Theorem 5.7]{St2} and also \cite[1.5]{JKSY2}. Note that $V_N$ (and not $V'_N$) is used in the definition of $B_{\Ga}$. By (\hl{1.3.5}{1.3.5}), we can get the spectrum from (\hl{1.4.3}{1.4.3}). Here it is not necessary to assume $f$ is {\it simplicial}.
\par\htt{R1.4a}{}\msn
{\bf Remark\,\,1.4a.} In the {\it non-simplicial\1} case, there is a counter-example to a conjecture of Steenbrink on the spectral pairs of Newton non-degenerate functions. For instance, if
$$f=x^2+y^2+xz+yz+z^4,$$
\par\nin \vskip-3mm\nin
then we have
$$q_{\si_1}=t+t^{3/2},\q\h{although}\q\Sp_f(t)=t^{3/2},$$
\par\nin where $\si_1$ is the unique {\it non-simplicial\1} 2-dimensional face of $\Gp(f)$, see \cite{Da2}. More precisely, the other $q_{\si}(t)$ vanish except for $q_{\emptyset}(t)=1$, and the combinatorial polynomials $r_{\si}(t)$ in (\hl{4}{4}) for these are given by $r_{\si_1}(t)=1$, $r_{\emptyset}(t)=-t$. This implies that the modified version of Steenbrink conjecture as in Theorem\,\,\hl{T2}{2} does not hold in the {\it non-simplicial\1} case.
\par\htt{R1.4b}{}\msn
{\bf Remark\,\,1.4b.} As a corollary of Theorem\,\,\hl{T2}{2}, the assertion (\hl{4}{4}) in the introduction holds also in the $f$ non-convenient case if $f$ is simplicial, and has an isolated singularity at 0.
For an example of a non-convenient function, we may consider for instance
$$f=x^ay\pl y^bz\pl xz^c\q(a,b,c\in\Z_{>0}).$$
\par\nin The intersection of the dual cone $\Si$ with the surface defined by $\msum_{i=1}^3\,v_i=1$ in $\R^3$ is given by
$$\setlength{\unitlength}{.4cm}
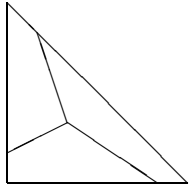
\begin{picture}(6,6)
\put(0,0){\line(0,1){6}}
\put(0,0){\line(1,0){6}}
\put(6,0){\line(-1,1){6}}
\put(2,2){\line(3,-2){3}}
\put(2,2){\line(-1,3){1}}
\put(2,2){\line(-2,-1){2}}
\end{picture}$$
\par\nin \sk
It does not seem necessarily easy to prove the above assertion without using Theorem\,\,\hl{T2}{2} by reducing to the convenient case (adding certain monomials to $f$ and applying \cite[Proposition A.2]{JKSY1} together with the finite determinacy of holomorphic functions with isolated singularities). For instance, we might have to show that the combinatorial polynomials $r_{\si}(t)$ do not change by passing from $f$ to $f\pl\msum_i\,x_i^{a_i}$ for $a_i\gg 0$ in the case $q_{\si}(t)\ne 0$. However, the combinatorics of non-convenient Newton polytopes can be rather complicated if $n\ges 4$.
\sk
One may consider for instance the following example with $n\eq 4$\,:
$$f=x^3y\pl y^3\pl x^2z\pl z^3\pl zw^3.$$
\par\nin This has an isolated singularity with Milnor number 37 according to Singular \cite{Sing}, and seems to be simplicial. Let $\si_1$ be the convex hull of $(3,1,0,0)$ and $(0,3,0,0)$. We have $q_{\si_1}(t)\ne 0$, and hence would have to calculate the combinatorial polynomial $r_{\si_1}(t)$ for $f$ and also for $g:=f\pl x^a\pl w^b$ ($a,b\gg 0$). Let $v^{(i)}\in\N^3$ be the lattice point corresponding to the $i$\1th monomial appearing in $g$ ($i\in[1,7]$). We would have to determine which subset of $\{v^{(i)}\}_{i\in[3,7]}$ gives the set of vertices of some face of $\Gg$ by taking the union with $\{v^{(1)},v^{(2)}\}$.
These seem to be $\{3,5\},\{5,7\},\{5\};\{3\},\{7\},\emptyset$, if the calculation is correct, where $v^{(i)}$ is simply denoted by $i$. Notice that there is a partition into two groups depending on whether $v^{(5)}$ is contained or not, and there is moreover a one-to-one corresponding between the two groups by adding $v^{(5)}$. This may be related closely to the invariance of combinatorial polynomials. Note also that, in the isolated singularity case, for each $i\in[1,n]$, there is $j\in[1,n]$ such that the coefficient of $x_i^{a_i}x_j$ in the Taylor expansion of $f$ does not vanish for some $a_i>0$ in general.
\bsn\msn\htt{S2}{}\par
\vbox{\centerline{\bf 2. Relation to the toric geometry}
\bsn
In this section we review some toric geometry, and prove the main theorems.}
\par\htt{2.1}{}\msn
{\bf 2.1.~Toric varieties.} Let $f:(\C^n,0)\to(\C,0)$ be a Newton non-degenerate holomorphic function with $\Gp(f)$ the Newton polytope. We have the {\it dual fan\1} $\Si$ in $(\R_{\ges 0})^n$, see for instance \cite[9.1]{Va1} (and Remark\,\,\hl{R2.1a}{2.1a} below). Let $\Xi$ be a smooth subdivision of $\Si$, see \cite[8.2]{Da1}, \cite{Od}, \cite{KKMS}. Here smooth means that any cone $\xi\in\Xi$ is simplicial and moreover the semi-group $\Z^n\cap\xi\setminus\{0\}$ is freely generated by the primitive elements in the 1-dimensional faces of $\xi$ so that $\C[\1\Z^n\cap\xi\1]$ is isomorphic to a polynomial ring. We have the toric varieties $\X$, $\Y$ associated to $\Si$, $\Xi$ respectively and also the proper morphisms
$$\Y\buildrel{\rho}\over\longrightarrow\X\buildrel{\pi}\over\longrightarrow X:=\C^n,$$
\par\nin inducing isomorphisms over $X^*:=(\C^*)^n$, where $X$ corresponds to the fan consisting of $(\R_{\ges 0})^n$ and its faces.
Note that the dual fan $\Si$ of $\Gp(f)$ contains every coordinate axis of $\R_{\ges 0}^n$. For any $(n{-}1)$-dimensional {\it compact\1} face $\si<\Gp(f)$, the corresponding 1-dimensional cone $\eta_{\si}$ is contained in $\R_{>0}^n\cup\{0\}$ by the {\it positivity\1} of the coefficients of a linear function defining $\si$. This implies that $\pi$ induces an isomorphism over $X\setminus\{0\}$ if $f$ is convenient.
\sk
In this paper we use the {\it associated analytic spaces\1} rather than algebraic varieties, since we have to deal with local systems defined on Zariski-open subsets of them. Let $U_0\subset X$ be a sufficiently small open neighborhood of $0\in X$ on which $f$ is defined. Set
$$\begin{array}{lll}X_0:=U_0\cap f^{-1}(0),&X^{(1)}:=\{x_1\cdots x_n=0\}\,\,\subset\,X,\\ \X_0:=\pi^{-1}(X_0),&\X^{(1)}:=\pi^{-1}(X^{(1)})\,\,\subset\,\,\X,\raise6mm\h{}\\ \Y_0:=\rho^{-1}(\X_0),&\Y^{(1)}:=\rho^{-1}(\X^{(1)})\,\,\subset\,\,\Y.\raise6mm\h{}\end{array}$$
\par\nin There are natural isomorphisms
$$\Y\setminus\Y^{(1)}\simto\X\setminus\X^{(1)}\simto X\setminus X^{(1)}=X^*:=(\C^*)^n.$$
\par\nin Note that $\Y$ is smooth, and $\Y^{(1)}\subset\Y$ is a divisor with normal crossings (since $\Xi$ is a smooth subdivision of $\Si$). Moreover $\Y_0\subset\rho^{-1}\pi^{-1}(U_0)$ is also a divisor with normal crossings (shrinking $U_0$ if necessary). Indeed, the intersection of the proper transform of $f^{-1}(0)$ with any positive-dimensional stratum $\Y_{\xi}$ of $\Y^{(1)}$ (which is a torus action orbit) is a smooth divisor on $\Y_{\xi}$, using the condition that $f$ is non-degenerate. Here the intersection is scheme-theoretic, since it is defined by restricting a local defining function $h$ in an ambient space so that its restriction to $\Y_{\xi}$ coincides (up to multiplication by a nowhere vanishing function) with the pull-back of $f_{\si}$ for $\si\les\Gf$ such that $\rho(\Y_{\xi})=\X_{\si}$ in the notation explained just below (and $\ddd h$ does not vanish after the restriction).
\sk
Put
$$\X_0^{\rm ex}:=\pi^{-1}(0)\subset\X_0,\q\Y_0^{\rm ex}:=\rho^{-1}(\X_0^{\rm ex})\subset\Y_0.$$
\par\nin Let $\X_0^{\rm pr}\subset\X_0$ be the proper transform of $X_0$ so that $\X_0=\X_0^{\rm ex}\cup\X_0^{\rm pr}$. By the definitions of the dual fan of $\Gp(f)$ and the associated toric varieties, there are stratifications
\htt{2.1.1}{}
$$\begin{array}{ll}\X^{(1)}=\msqcup_{\si<\Gp(f)}\,\X_{\si},\q&\X_0^{\rm ex}=\msqcup_{\si\les\Gf}\,\X_{\si},\\ \Y^{(1)}=\msqcup_{\xi\1\in\1\Xi}\,\Y_{\xi},&\Y_0^{\rm ex}=\msqcup_{\xi\1\in\1\Xi'}\,\Y_{\xi},\raise6mm\h{}\end{array}
\leqno(2.1.1)$$
\par\nin where $\Xi'\subset\Xi$ consists of cones which are not contained in any coordinate hyperplane of $\R^n$, see also Remark\,\,\hl{R2.1b}{2.1b} below. Note that the above stratification of the exceptional divisor $\X_0^{\rm ex}\eq\pi^{-1}(0)$ is \h{\it strictly\1} \h{\it finer\1} than the usual one used for the construction of motivic nearby fibers in \cite[3.3]{DL} (for instance, if $f\eq\msum_{i=1}^n\,x_i^d$ ($d\ges 2$), we have $\pi^{-1}(0)\eq\PP^{n-1}$, which is smooth). The numbers $k(\si)$ in Theorem\,\,\hl{T1}{1} is closely related to this difference. Note that the number of irreducible components of $\Y_0$ passing through a point of $\Y_{\xi}^{\circ}$ is $d_{\xi}-\kc_{\xi}$ in the notation of (\hl{2.3.1}{2.3.1}) below.
\sk
The stratifications of $\X^{(1)}$, $\X_0^{\rm ex}$ can be indexed also by $\Si$, $\Si'$ respectively, where $\X_{\si}\eq\X_{\eta_{\si}}$ in the notation of Remark\,\,\hl{R2.1a}{2.1a} below (and $\Si'\subset\Si$ is defined in the same way as $\Xi'\subset\Xi$).
Note that $\X_{\si}\cong(\C^*)^{d_{\si}}$, $\Y_{\xi}\cong(\C^*)^{n-d_{\xi}}$ for $\si<\Gp(f)$, $\xi\in\Si$, and they are {\it orbits\1} of the natural torus actions on $\X$, $\Y$, see Remarks\,\,\hl{R2.1a}{2.1a}--\hl{R2.1c}{c} below. We have $\si\ges\tau$ (or equivalently $\eta_{\si}\les\eta_{\tau}$) if and only if the closure $\Xo_{\si}$ contains $\X_{\tau}$ (and similarly for $\Y_{\xi}$).
\par\htt{R2.1a}{}\msn
{\bf Remark\,\,2.1a.} For $v\in\R_{\ges 0}^n$, set
$$\Gp(f)_v:=\bl\{u\in\Gp(f)\,\big|\,\langle v,u\rangle=\min\langle v,\Gp(f)\rangle\br\}$$
\par\nin where $\langle v,u\rangle$ denotes the natural pairing of $v,u\in\R^n$.
Then $\Gp(f)_v$ is a face of $\Gp(f)$. For a face $\si\les\Gp(f)$, set
$$\eta_{\si}:=\bl\{v\in\R_{\ges 0}^n\,\big|\,\Gp(f)_v\supset\si\br\}\,\subset\,\R_{\ges 0}^n.$$
\par\nin This defines the dual fan $\Si$ to $\Gp(f)$, see for instance \cite[9.1]{Va1}. (It is also called the normal fan, where min is replaced by max usually; so the dual cones have {\it opposite\1} directions.)
\sk
By this construction, each $(n{-1})$-dimensional (not necessarily compact) face $\si$ of $\Gp(f)$ corresponds to a 1-dimensional dual cone $\eta_{\si}\in\Si$, and then to the closure $\Xo_{\si}$ of $\X_{\si}$ in $\X$, which is the divisor corresponding to $\eta_{\si}$.
\sk
For a (not necessarily compact) face $\si<\Gp(f)$, let $\tau_1,\dots,\tau_m$ be the $(n{-1})$-dimensional faces of $\Gp(f)$ containing $\si$ with $\eta_{\tau_i}\in\Si$ the 1-dimensional dual cone to $\tau_i$ ($i\in[1,m]$), where $m\ges n-\dim\si$. Then $\si$ corresponds to the cone $\eta_{\si}\in\Si$ spanned by $\eta_{\tau_1},\dots,\eta_{\tau_m}$, and then to the intersection $\Xo_{\tau_1}\cap\dots\cap\Xo_{\tau_m}$ in $\X$.
\sk
The latter assertion means that the combinatorial data of the set of {\it intersections of divisors\1} $\Xo_{\si}$ with $d_{\si}=n{-}1$ are the same as those of the set of (not necessarily compact) faces of $\Gp(f)$. Here the intersection $\Xo_{\tau_1}\cap\dots\cap\Xo_{\tau_m}$ with $d_{\tau_i}=n{-}1$ is {\it non-compact\1} if and only if $\si=\mcap_i\,\tau_i$ is a {\it non-compact\1} face of $\Gp(f)$, or equivalently, all the $\eta_{\tau_i}$ is contained in one coordinate hyperplane of $\R^n$ (since $\eta_{\si}$ does not intersect $\R_{>0}^n$, and is spanned by the $\eta_{\tau_i}$).
\par\htt{R2.1b}{}\msn
{\bf Remark\,\,2.1b.} For $I\subset\{1,\dots,n\}$, set
$$\aligned\XXo_I&:=\mcap_{i\in I}\{x_i\eq 0\}\subset X,\\
X_I&:=\XXo_I\setminus\mcup_{i\notin I}\,\XXo_{I\cup\{i\}},\\
\ga_I&:=\begin{cases}0&\h{if}\q U_0\cap X_I\subset f^{-1}(0),\\
1&\h{otherwise.}\end{cases}\endaligned$$
\par\nin The $X_I$ are the strata of the natural stratification of $X=\C^n$, which is compatible with the morphism of toric varieties $\pi:\X\to X$. If there is a face $\si<\Gp(f)$ such that $\pi(\X_{\si})=X_I$ and $d_{\si}>d_I$, then we have $\ga_I=0$, where $d_I:=n-|I|$. This is verified inductively using the projections as in \cite[A.1.4--5]{JKSY1}.
Note that the first condition $\pi(\X_{\si})=X_I$ is equivalent to that $I\subset\{1,\dots,n\}$ is the maximal subset such that
$$\eta_{\si}\subset\R_{\ges 0,I}\,\bl(:=\mcap_{i\in I}\,\{v_i=0\}\subset\R_{\ges 0}^n\br),$$
\par\nin which is denoted by $I_{\si}$. (This is similar to the assertion for $\rho$ at the end of Remark\,\,\hl{R2.1d}{2.1d} below.) We have
$$k(\si)=n-|I_{\si}\setminus J_f|.$$
\par\nin $$I_{\si}=1\q\h{if}\q I_{\tau}=1,\,\,\si<\tau.$$
\par\nin \par\htt{R2.1c}{}\msn
{\bf Remark\,\,2.1c.} Let $U_{\si}\subset\X$ be the affine open subset corresponding to $\si<\Gp(f)$, or equivalently to $\eta_{\si}\in\Si$. We have
\htt{2.1.2}{}
$$\aligned U_{\si}&={\rm Spec}\,\C[\Z^n\cap\eta_{\si}^{\vee}]^{\rm an}=\msqcup_{\si\les\tau\les\Gp(f)}\,\X_{\tau}\,,\\ \X_{\si}&={\rm Spec}\,\C[\Z^n\cap\eta_{\si}^{\perp}]^{\rm an},\endaligned
\leqno(2.1.2)$$
\par\nin (where ${\rm Spec}\,\A^{\rm an}$ is the analytic space associated to ${\rm Spec}\,\A$) with
$$\aligned\eta_{\si}^{\vee}:=\bl\{u\in\R^n\mid\langle v,u\rangle\ges 0\,\,(\forall\,v\in\eta_{\si})\br\},\\ \eta_{\si}^{\perp}:=\bl\{u\in\R^n\mid\langle v,u\rangle=0\,\,(\forall\,v\in\eta_{\si})\br\}.\endaligned$$
\par\nin We have $\eta_{\si}^{\vee}=\mcap_i\,\eta_{\tau_i}^{\vee}$ with $\eta_{\tau_i}$ as in Remark\,\,\hl{R2.1a}{2.1a}, and
$$U_{\si}\supset\X_{\si_0}:=X^*\cong(\C^*)^n\q\h{with}\q\si_0:=\Gp(f).$$
\par\nin \sk
Notice that $\C[\Z^n\cap\eta_{\si}^{\perp}]$ is a {\it subquotient\1} of $\C[\Z^n]$ by the injective and surjective morphisms
$$\C[\Z^n]\hookleftarrow\C[\Z^n\cap\eta_{\si}^{\vee}]\onto\C[\Z^n\cap\eta_{\si}^{\perp}],$$
\par\nin corresponding to the open and closed immersions of analytic spaces
$$(\C^*)^n\into U_{\si}\hookleftarrow\X_{\si}.$$
\par\nin \par\htt{R2.1d}{}\msn
{\bf Remark\,\,2.1d.} Let $V_{\xi}\subset\Y$ be the affine open subset corresponding to $\xi\in\Xi$. We have
\htt{2.1.3}{}
$$\aligned V_{\xi}&={\rm Spec}\,\C[\Z^n\cap\xi^{\vee}]^{\rm an}=\msqcup_{\xi'\les\xi}\,\Y_{\xi'}\,,\\ \Y_{\xi}&={\rm Spec}\,\C[\Z^n\cap\xi^{\perp}]^{\rm an}\endaligned
\leqno(2.1.3)$$
\par\nin together with the open and closed immersions
$$(\C^*)^n\into V_{\xi}\hookleftarrow\Y_{\xi}.$$
\par\nin The above disjoint union is taken over any $\xi'\les\xi$ including the case $\xi'=\emptyset$ with $\Y_{\emptyset}=(\C^*)^n$.
\sk
Note that $\rho$ induces the natural morphism $V_{\xi}\to U_{\eta}$ if $\xi\subset\eta$, or equivalently, if $\xi^{\vee}\supset\eta^{\vee}$. The morphism $\rho$ induces the smooth surjective morphism $\Y_{\xi}\to\X_{\eta}$ (that is $\rho(\Y_{\xi})=\X_{\eta}$) if and only if $\eta\in\Si$ is the minimal dimensional cone containing $\xi\in\Xi$ (which will be denoted by $\eta_{\xi}$). These morphisms are compatible with (\hl{2.1.2}{2.1.2}--\hl{2.1.3}{3}).
\par\htt{R2.1e}{}\msn
{\bf Remark\,\,2.1e.} There are natural isomorphisms of analytic spaces:
\htt{2.1.4}{}
$$\X_0^{\rm ex}=({\rm Proj}\,A_{\Ga})^{\rm an},\q\Xo_{\si}=({\rm Proj}\,A_{\si})^{\rm an}\,\,\,\,(\si\les\Gf),
\leqno(2.1.4)$$
\par\nin in the notation of \hl{1.4}{1.4}, see \cite{Ha} for ${\rm Proj}$. It is easy to see the first isomorphism locally on the complement of the divisor on $({\rm Proj}\,A_{\Ga})^{\rm an}$ defined by $x^{\nu}$ for $\si=\{\nu\}\in\CFf^0$, since the latter is naturally isomorphic to the closed subvariety of $U_{\si}$ in (\hl{2.1.2}{2.1.2}) defined by the ideal generated by $x^{\nu}$ with $\nu\in\Z^n$ contained the {\it interior\1} of $\eta_{\si}^{\vee}$. Note that for $[x^{\nu}]$, $[x^{\nu'}]\in A_{\Ga}$, we have
$$[x^{\nu}]{\cdot}[x^{\nu'}]=0\,\,\,\,\h{unless}\,\,\,\,\nu,\nu\1'\in C(\tau)\,\,\,\,\h{for some}\,\,\,\,\tau\les\Gf.$$
\par\nin A similar argument shows the second isomorphism of (\hl{2.1.4}{2.1.4}).
(Note that $A_{\si}$ is a quotient ring of $A_{\Ga}$, and ${\rm Proj}\,A_{\si}$ is an irreducible component of ${\rm Proj}\,A_{\Ga}$ if $d_{\si}=n{-}1$.) The latter implies that $\Xo_{\si}\setminus\X_0^{\rm pr}$ is an {\it affine\1} variety, since its complement $\Xo_{\si}^{\rm pr}:=\Xo_{\si}\cap\X_0^{\rm pr}\subset\Xo_{\si}$ is defined by $f_{\si}\in A_{\si}$.
\par\htt{R2.1f}{}\msn
{\bf Remark\,\,2.1f.} The motivic nearby fiber defined as in \hl{A.2}{A.2} below can be obtained up to multiplication by $(1{-}\Lb)^{k(\si)-d(\si)}$ by restricting the following over $X_{\si}^{\circ}\subset\X_0^{\rm ex}$\,:
\htt{2.1.5}{}
$${\rm Spec}\bl(A_{\Ga}/A_{\Ga}(f_{\Gf}{-}1)\br)^{\rm an}\q\h{over}\q\X_0^{\rm ex}.
\leqno(2.1.5)$$
\par\nin Here we use the identification
$$\Gr_{V'_N}^1f\cong f_{\Gf}.$$
\par\nin Taking the base changes of (\hl{2.1.5}{2.1.5}), we get
\htt{2.1.6}{}
$$\aligned{\rm Spec}\bl(A_{\si}/A_{\si}(f_{\si}{-}1)\br)^{\rm an}\q\h{over}\q\Xo_{\si},\\{\rm Spec}\bl(\Ah_{\si}/\Ah_{\si}(f_{\si}{-}1)\br)^{\rm an}\q\h{over}\q X_{\si},\endaligned
\leqno(2.1.6)$$
\par\nin where $\Ah_{\si}$ is the subring of $\C[\Z^n]$ generated by $A_{\si}=\C[\Z^n\cap C(\si)]$ and $\C[\Z^n\cap\eta_{\si}^{\perp}]$ (with $C(\si)$ the cone of $\si$). The last term of (\hl{2.1.6}{2.1.6}) can be identified with $\Xt_{\si}^{\circ}$ defined in (\hl{2.2.2}{2.2.2}) below, since $\X_{\si}^{\rm pr}\subset\X_{\si}$ is defined by $h_{\si}:=f_{\si}/x^{\nu'_{\si}}$ with $\nu^{\1\prime}_{\si}$ explained just below, see also the beginning of Section\,\,\hl{2.2}{2.2} below.
\sk
Notice that $({\rm Spec}\,\Ah_{\si})^{\rm an}$ is a line bundle over $X_{\si}$.
We get a $\C^*$-bundle if $\Ah_{\si}$ is replaced by $\C[\Z^n\cap V(\si)]$ with $V(\si)\subset\R^n$ the vector subspace spanned by $\si$.
Their {\it trivializations\1} are given by choosing a monomial $x^{\nu_{\si}}$ such that
\htt{2.1.7}{}
$$\nu^{\1\prime}_{\si}:=\de_{\si}\1\nu_{\si}\,\in\,A(\si)=\si\pl\eta_{\si}^{\perp}.
\leqno(2.1.7)$$
\par\nin Recall that $A(\si)$ is the smallest affine space containing $\si$, see Convention at the end of the introduction. We can identify $x^{\nu_{\si}}$ with a coordinate of the line bundle, see also a remark after (\hl{2.2.2}{2.2.2}). This kind of argument seems to be needed to relate motivic nearby fibers of Denef and Loeser to certain constructions in some papers quoted by \cite{Sta} as relevant results.
\par\htt{R2.1g}{}\msn
{\bf Remark\,\,2.1g.} For a holomorphic function $g$ on a complex manifold $Y$ with $Y_0:=g^{-1}(0)$ a divisor with simple normal crossings, it is well known that there is a {\it normal\1} analytic space $\Y_0$ over $Y_0$ such that the motivic nearby fiber is given up to multiplication by a power of $1\mi\Lb$ by its restriction over each stratum of the stratification associated with the normal crossing divisor. This is obtained by taking the normalization of the base change of $g:Y\to\C$ by a ramified covering $\C\ni\ttt\mapsto t=\ttt\1^m\in\C$ for a sufficiently divisible positive integer $m$ (more precisely, $m$ is divisible by the multiplicity of any irreducible component of $Y_0$), see, for instance, an argument before \cite[Lemma 5.3]{Loo} and also \cite{St1}. Note that $\Y_0$ is {\it algebraic\1} if $Y$, $g$ are.
\sk
It does not seem trivial to show that the above motivic nearby fiber is isomorphic to the one constructed in \hl{2.2}{2.2} below (since the normalization is involved). At the level of nearby cycle sheaves, however, it is rather easy to show the isomorphism, since the universal covering of $\C^*$ (used for the definition of nearby cycles sheaves) factors through the above finite morphism (restricted over $\C^*$). In our case, we can then verify the isomorphism at the level of finite unramified covering space over each $\X_{\si}^{\circ}$ by applying Remark\,\,\hl{R2.1h}{2.1h} as well as Remarks\,\,\hl{R2.2a}{2.2a}--\hl{R2.2b}{b} below (using also the faithfulness of the permutation representation of the symmetric group associated with a finite unramified covering). Note that $T_s$ can be identified with the action of a generator of the cyclic covering transformation group, and the local monodromy around a divisor at infinity is a power of it by Remark\,\,\hl{R2.2d}{2.2d} below.
\sk
The variety in (\hl{2.1.5}{2.1.5}) may be normal. However, it is {\it not} defined over $\X_0^{\rm ex}\cap\X_0^{\rm pr}$, and we need an argument as in Remark\,\,\hl{R2.2b}{2.2b} below.
\par\htt{R2.1h}{}\msn
{\bf Remark\,\,2.1h.} Unramified finite covering spaces of smooth (or more generally, normal) complex algebraic varieties can be determined uniquely by their underlying topological (or equivalently, analytic) spaces. This means that any topological (or analytic) unramified finite covering space of a smooth complex variety has a unique structure as algebraic variety, see \cite{GR}, \cite[XII, Cor.~5.2]{Gr2}, \cite[App.~B, 3.2]{Ha}, etc. Here an unramified finite covering means a finite \'etale morphism in the algebraic case, and a locally isomorphic finite morphism in the topological (or analytic) case because of the difference in topology. (Note that a finite morphism is a {\it proper\1} morphism with finite fibers.) The above assertion is sometimes called the generalized Riemann existence theorem.
\par\htt{2.2}{}\msn
{\bf 2.2.~Construction of nearby cycle sheaves.} For $\si<\Gp(f)$, set
$$\X_{\si}^{\rm pr}:=\X_{\si}\cap\X_0^{\rm pr},\q\q\X_{\si}^{\circ}:=\X_{\si}\setminus\X_{\si}^{\rm pr},$$
\par\nin with $i_{\si}:\X_{\si}^{\rm pr}\into\X_{\si}$, $j_{\si}:\X_{\si}^{\circ}\into\X_{\si}$ the inclusions.
Note that $\X_{\si}^{\rm pr}\subset\X_{\si}\cong(\C^*)^{d_{\si}}$ is smooth for $\si\les\Gf$, since it is defined {\it essentially\1} by $f_{\si}$, more precisely, by $h_{\si}:=f_{\si}/x^{\nu'_{\si}}$ with $\nu^{\1\prime}_{\si}$ as in the end of Remark\,\,\hl{R2.1f}{2.1f}.
\sk
We have the equality of Weil divisors on the normal variety $\X$\,:
\htt{2.2.1}{}
$${\rm div}\,\pi^*f=\msum_{d_{\tau}=n-1}\,\de_{\tau}\1\Xo_{\tau}+\X_0^{\rm pr},
\leqno(2.2.1)$$
\par\nin with $d_{\tau}:=\dim\tau$. This follows from the definition of $\de_{\tau}$, see (\hl{1}{1}) in the introduction. We see that the first term of the right-hand side is locally principal (hence so is the last term), considering the pull-back of any functions $g$ with $\Gp(g)\subset\Gp(f)$ (see also \cite{Da1}, \cite{Od} for the general theory of equivariant invertible sheaves on toric varieties).
\sk
Replacing the last condition with $\Gp(g)\subset\tfrac{\!\!1}{\de_{\si}}\1\Gp(f)$, we see that the restriction of
$$D_{\si}:=\msum_{d_{\tau}=n-1,\tau\ges\si}\,(\de_{\tau}/\de_{\si})\1\Xo_{\tau}$$
\par\nin to a Zariski-open neighborhood of $\X_{\si}$ is locally principal. Note that $\de_{\tau}$ is divisible by $\de_{\si}$ if $\tau\ges\si$. Consider the invertible sheaf $\Lc_{\si}$ on $\X_{\si}$ defined by
$$\Lc_{\si}:=\OO_{\X}(D_{\si})|_{\X_{\si}}.$$
\par\nin The equality of divisors (\hl{2.2.1}{2.2.1}) then gives an isomorphism of invertible sheaves on $\X_{\si}$\,:
$$\Lc_{\si}^{\otimes\de_{\si}}\cong\OO_{X_{\si}}(-\X_{\si}^{\rm pr}).$$
\par\nin We thus get a cyclic covering $p_{\si}:\Xt_{\si}\to\X_{\si}$ defined by
$$\Xt_{\si}:=\Sc pec_{\,\X_{\si}}\bl(\mopl_{k=0}^{\de_{\si}-1}\,\Lc_{\si}^{\otimes k}\br)^{\rm an},$$
\par\nin see also \cite[Appendix]{wt0}. This is totally ramified along $\X_{\si}^{\rm pr}$, and $\Xt_{\si}^{\circ}:=p_{\si}^{-1}(\X_{\si}^{\circ})$ is a finite unramified covering of $\X_{\si}^{\circ}$. Indeed, $\Xt_{\si}$ is defined by
\htt{2.2.2}{}
$$z^{\de_{\si}}=h_{\si}\q\h{in}\q\C\times\X_{\si}.
\leqno(2.2.2)$$
\par\nin Here $z$ is the coordinate of $\C$, and $h_{\si}$ is a defining function of $\X_{\si}^{\rm pr}$ in $\X_{\si}\cong(\C^*)^{d_{\si}}$ explained in the beginning of this section. This construction is essentially the same as motivic nearby fibers of Denef and Loeser (see \cite[3.3]{DL}) except that $z$ in (\hl{2.2.2}{2.2.2}) is replaced by $z^{-1}$. This replacement is reasonable when one defines the motivic nearby fibers as in Remark\,\,\hl{R2.1f}{2.1f}, since $z^{-1}$ is essentially identified with $x^{\nu_{\si}}$ in (\hl{2.1.7}{2.1.7}).
\sk
We define a $T_s$-action on $\Xt_{\si}$ over $\X_{\si}$ choosing a generator of the covering transformation group $\Z/\de_{\si}\Z$, which acts on the above coordinate $z$ by the multiplication by $\exp(2\pi\sqrt{-1}/\de_{\si})$.
We then get a $\Q$-local system of rank $\de_{\si}$ on $\X_{\si}^{\circ}$\,:
$$L_{\si}:=(p_{\si})_*\Q_{\Xt_{\si}}|_{\X_{\si}^{\circ}},$$
\par\nin which is endowed with an action of $T_s$ of order $\de_{\si}$.
This is viewed as a variation of Hodge structure of type $(0,0)$ on $\X_{\si}^{\circ}$.
\sk
Setting
$$L_{\si,\C,\la}:={\rm Ker}(T_s\mi\la)\subset L_{\si,\C}:=\C\otimes_{\Q}L_{\si},$$
\par\nin we have the direct sum decomposition
\htt{2.2.3}{}
$$L_{\si,\C}=\mopl_{\la\in\mu_{\de_{\si}}}\,L_{\si,\C,\la},
\leqno(2.2.3)$$
\par\nin with ${\rm rank}\,L_{\si,\C,\la}=1$ if $\la\in\mu_{\de_{\si}}$ (and 0 otherwise). Here $\mu_m:=\{\la\in\C\mid\la^m=1\}$ ($m\ges 2$).
\sk
With $\Q$-coefficients, we have the decomposition
\htt{2.2.4}{}
$$\aligned&\q\q L_{\si}=L_{\si,1}\oplus L_{\si,\1\ne1}\q\q\h{with}\\&L_{\si,1}:={\rm Ker}(T_s\mi{\rm id})=\Q_{\X_{\si}^{\circ}}\subset L_{\si},\endaligned
\leqno(2.2.4)$$
\par\nin in a compatible way with the action of $T_s$ (since $p_{\si}$ is totally ramified along a smooth hypersurface $\X_{\si}^{\rm pr}\subset\X_{\si}$).
\par\htt{R2.2a}{}\msn
{\bf Remark\,\,2.2a.} The local systems $L_{\si,\C,\la}$ are uniquely determined by {\it local monodromies\1} around divisors at infinity of a smooth compactification of $\X_{\si}^{\circ}$ such that the complement of $\X_{\si}^{\circ}$ is a divisor with normal crossings. Indeed, if there are two local systems $L_1,L_2$ of rank~1 on $\X_{\si}^{\circ}$ with the same local monodromies at infinity, we can consider the tensor product $L_1\otimes L_2^{\vee}$. This has trivial local monodromies at infinity, and can be extended to a local system on the smooth compactification. But it must be a trivial local system, since $\X_{\si}^{\circ}$ has a simply connected smooth partial compactification $\C^{d_{\si}}$.
\par\htt{R2.2b}{}\msn
{\bf Remark\,\,2.2b.} For $\si\les\Gf$, there are distinguished triangles
\htt{2.2.5}{}
$$(j_{\si})_!L_{\si}\to\RR(j_{\si})_*L_{\si}\to(i_{\si})_*i_{\si}^*\RR(j_{\si})_*L_{\si}\buildrel{+1}\over\to,
\leqno(2.2.5)$$
\par\nin together with isomorphisms of $\Q$-local systems on $\X_{\si}^{\rm pr}:$
\htt{2.2.6}{}
$$i_{\si}^*R^k(j_{\si})_*L_{\si}=i_{\si}^*R^k(j_{\si})_*\Q_{\X_{\si}^{\circ}}=\begin{cases}\Q_{\X_{\si}^{\rm pr}}&(k=0),\\ \Q_{\X_{\si}^{\rm pr}}(-1)&(k=1),\\ \,0&(k\ne0,1),\end{cases}
\leqno(2.2.6)$$
\par\nin since $p_{\si}$ is totally ramified along a smooth hypersurface $\X_{\si}^{\rm pr}\subset\X_{\si}$. These are needed to show the compatibility of (\hl{2}{2}) with a formulation using motivic nearby fibers in Theorem\,\,\hl{TA.2}{A.2} below. Theorem\,\,\hl{T1}{1} is also compatible with \cite[1.3]{BS2} in the case of homogeneous polynomials with an isolated singularity at 0. These are closely related to \cite[Theorem 4.2]{BS1}.
\msn
\sk
We have no problems as above for the {\it non-unipotent\1} monodromy part $L_{\si,\1\ne1}$, since there are canonical isomorphisms
\htt{2.2.7}{}
$$(j_{\si})_!L_{\si,\1\ne1}=(j_{\si})_*L_{\si,\1\ne1}=\RR(j_{\si})_*L_{\si,\1\ne1},
\leqno(2.2.7)$$
\par\nin and these give the {\it intermediate direct image\1} of $L_{\si,\1\ne1}$ in \cite{BBD} (up to a shift of complex).
\par\htt{R2.2c}{}\msn
{\bf Remark\,\,2.2c.} For any local system $L$ of rank 1 on a torus $T$, we have
$$H^j(T,L)=0\q(\forall\,j\in\Z),$$
\par\nin if and only if $L$ is not a trivial local system (that is, some local monodromies around a divisor at infinity is non-trivial). In the trivial case, we have in the Grothendieck group of mixed $\Q$-Hodge structures
$$\chi_c(T,\Q)=(\theta\mi1)^{\dim T}.$$
\par\nin \par\htt{R2.2d}{}\msn
{\bf Remark\,\,2.2d.} Let $g$ be a function on a complex manifold $Y_1$. Set $f=gz^m$ on $Y=Y_1\times\De$ with $z$ the coordinate of $\De$ and $m\in\Z_{>0}$. Then the monodromy around $z=0$ of the nearby cycle sheaves of $f$ on $Y_1\times\De^*$ is given by $T^{-m}$ with $T$ the monodromy of the nearby cycles (which is the inverse of the Milnor monodromy, see \cite{DS}).
This is compatible with a calculation of nearby cycles in the normal crossing case \cite[3.3.1]{mhm} (comparing $\nu$ and $\nu\1'\pl\nu_0m$.)
\par\htt{R2.2e}{}\msn
{\bf Remark\,\,2.2e.} For $\xi\in\Xi'$, let $\xi_i$ be its 1-dimensional faces ($0\in[1,d_{\xi}]$). Let $\si\les\Gf$ be the (unique) face such that its corresponding $\eta_{\si}\in\Si'$ coincides with $\eta_{\xi}$, that is, the minimal dimensional cone containing $\xi$, see Remark\,\,\hl{R2.1d}{2.1d}. Set
$$\de_{\si,\1\xi}:=\de_{\si\pl\xi^{\perp}},$$
\par\nin and similarly for $\de_{\si,\1\xi_i}$, see Convention at the end of the introduction and (\hl{2.1.7}{2.1.7}). (Note that $\si\pl\xi^{\perp}\supset\si\pl\eta_{\si}^{\perp}=A(\si)$.) We then get the equality
\htt{2.2.8}{}
$$\de_{\si,\1\xi}={\rm GCD}\bl(\de_{\si,\1\xi_i}\,\big|\,\,i\in[1,d_{\xi}]\br),
\leqno(2.2.8)$$
\par\nin since $\Y_0^{\rm ex}\subset\Y$ is a divisor with normal crossings, and the {\it multiplicity\1} of the pull-back of $f$ along the divisor $\Y_{\xi_i}$ is given by $\de_{\si,\1\xi_i}$ using (\hl{2.1.3}{2.1.3}). (Note that $\xi^{\perp}=\mcap_i\,\xi_i^{\perp}$.) In particular, $\de_{\si,\1\xi_i}/\de_{\si,\1\xi}\in\N$. We have $\de_{\si,\1\xi_i}=0$ if $\xi_i$ is a coordinate axis of $\R_{\ges 0}^n$. (We define GCD via the sum of of the corresponding ideals of $\Z$ so that $\Z\1\de_{\si,\1\xi}=\msum_i\,\Z\1\de_{\si,\1\xi_i}$, using a very classical theorem assuring that any ideal of $\Z$ is generated by one element.)
\par\htt{2.3}{}\msn
{\bf 2.3.~Proof of Theorem\,\,\hl{T1}{1}.} By the same argument as in \hl{2.2}{2.2}, we can define for $\xi\in\Xi'$ the ramified covering spaces
$$p_{\xi}:\Yt_{\xi}\to\Y_{\xi},$$
\par\nin together with the $\Q$-local systems $L'_{\xi}$ on
$$\Y_{\xi}^{\circ}:=\Y_{\xi}\setminus\Y_{\xi}^{\rm pr}\q\h{with}\q\Y_{\xi}^{\rm pr}:=\Y_{\xi}\cap\Y_0^{\rm pr}\,,$$
\par\nin where $\Y_0^{\rm pr}\subset\Y_0$ is the proper transform of $X_0$ in the notation of \hl{2.1}{2.1}. This essentially coincides with the construction of motivic nearby fibers of Denef and Loeser \cite[3.5.3]{DL}.
\sk
Let $\kc_{\xi}$ be the number of unit vectors $e_i\in\R^n$ contained in $\xi\in\Xi'$ with $i\notin J_f$ (similarly for $\kc_{\eta}$ with $\eta\in\Si'$), where $J_f$ is as in the introduction. By an argument similar to \cite[5.1 and 5.3]{MSS} and using \cite[Theorem 4.2]{BS1} (see also Remark\,\,\hl{R2.2b}{2.2b}), we get the equality in $K_0^{T_s}({\rm MHS})$
\htt{2.3.1}{}
$$[\chi(\Ff,\Q),T_s]=\msum_{\xi\1\in\1\Xi'}\,(1\mi\theta)^{d_{\xi}-\kc_{\xi}-1}\1\psi'_{\xi},
\leqno(2.3.1)$$
\par\nin with
$$\psi'_{\xi}:=\chi_c\bl(\Y_{\xi},\RR(j'_{\xi})_*(L'_{\xi},T_s)\br).$$
\par\nin Here $j'_{\xi}:\Y_{\xi}^{\circ}\into\Y_{\xi}$ is the inclusion, and $\chi_c$ denotes the Euler characteristic with compact supports. Note that $d_{\xi}\mi \kc_{\xi}$ is the number of irreducible components of $\Y_0$ passing through each point of $\Y_{\xi}^{\circ}$. In the notation of Remark\,\,\hl{R2.1a}{2.1a}, we have the equality
\htt{2.3.2}{}
$$k(\si)\mi d(\si)=d_{\eta_{\si}}\mi \kc_{\eta_{\si}}\mi 1,
\leqno(2.3.2)$$
\par\nin since
$$k(\si)\pl\kc_{\eta_{\si}}=n=d_{\si}\pl d_{\eta_{\si}}.$$
\par\nin \sk
We have
\htt{2.3.3}{}
$${\rm rank}\,L'_{\xi}=\de_{\si,\1\xi},
\leqno(2.3.3)$$
\par\nin that is, the degree of the ramified covering $p_{\xi}$ is $\de_{\si,\1\xi}$, see Remark\,\,\hl{R2.2e}{2.2e}. Here $\si\les\Gf$ is the face such that the corresponding $\eta_{\si}\in\Si'$ coincides with $\eta_{\xi}$ (the minimal dimensional cone containing $\xi$, see Remark\,\,\hl{R2.1d}{2.1d}), and $\de_{\si,\1\xi}$ is as in Remark\,\,\hl{R2.2e}{2.2e}. Then $\rho$ induces the smooth surjective morphism
$$\rho_{\xi}:\Y_{\xi}\onto\X_{\si},$$
\par\nin whose fibers are {\it tori}. (This is induced by (\hl{2.1.2}{2.1.2}--\hl{2.1.3}{3}).) We have
$$\de'_{\xi}:=\de_{\si,\1\xi}/\de_{\si}\in\Z,$$
\par\nin using the definition of $\de_{\si}$ as in (\hl{1}{1}), since there is an inclusion of affine spaces
$$\si\pl\xi^{\perp}\supset\si\pl\eta_{\si}^{\perp}.$$
\par\nin \sk
By (\hl{2.2.2}{2.2.2}) the pull-back of the ramified covering space $\Xt_{\si}\to\X_{\si}$ by $\rho_{\xi}$ is identified with a {\it quotient\1} covering space $\Ga'_{\xi}$ of $p_{\xi}:\Yt_{\xi}\to\Y_{\xi}$ defined by the (unique) cyclic subgroup of order $\de'_{\xi}$ of the covering transformation group $\Ga_{\xi}$ of $p_{\xi}$. Note that
$$\Ga'_{\xi}\cong\Z/\Z\de'_{\xi},\q\Ga_{\xi}\cong\Z/\Z\de_{\si,\1\xi},\q\Ga_{\xi}/\Ga'_{\xi}\cong\Z/\Z\de_{\si}.$$
\par\nin Then $\rho_{\xi}^*(L_{\si})$ is identified with the $\Ga'_{\xi}$-{\it invariant\1} part of $L'_{\xi}$, and we get the isomorphisms
\htt{2.3.4}{}
$$L'_{\xi,\C,\la}=\rho_{\xi}^*(L_{\si,\C,\la})\q\h{for any}\,\,\,\la\in\mu_{\de_{\si}}.
\leqno(2.3.4)$$
\par\nin \sk
Using these, we can show the following equality in $K_0^{T_s}({\rm MHS})$\,:
\htt{2.3.5}{}
$$\chi_c\bl(\Y_{\xi},\RR(j'_{\xi})_*(L'_{\xi},T_s)\br)=(\theta\mi1)^{r(\xi)}\1\chi_c\bl(\X_{\si},\RR(j_{\si})_*(L_{\si},T_s)\br),
\leqno(2.3.5)$$
\par\nin where $r(\xi):=\dim\Y_{\xi}\mi\dim\X_{\si}\,(=n\mi d_{\xi}\mi d_{\si})$. By Remark\,\,\hl{R2.2c}{2.2c}, it is enough to show that some local monodromy at infinity of the restriction of $L'_{\xi,\C,\la}$ to a general fiber of the projection $\rho_{\xi}^{\circ}:\Y_{\xi}^{\circ}\to\X_{\si}^{\circ}$ (which is a torus) is non-trivial if $\la\in\mu_{\de_{\si,\xi}}\setminus\mu_{\de_{\si}}$. But this can be verified by using Remarks\,\,\hl{R2.2d}{2.2d}--\hl{R2.2e}{e} as follows.
\sk
Let $\xi'\in\Xi'$ be a maximal dimensional cone containing $\xi$ and contained in $\eta_{\si}$, where $\si$ is as in Remark\,\,\hl{R2.2e}{2.2e}. Let $\xi'_j$ ($j\in[1,d_{\xi'}{-}d_{\xi}]$) be the 1-dimensional cones of $\xi'$ {\it not\1} contained in $\xi$. They correspond to normal crossing divisors $\Yo_{\xi'_j}$ meeting at $\Y_{\xi'}$. Their intersections with the closure $\Yo_{\xi}$ of $\Y_{\xi}$ are divisors at infinity of $\Y_{\xi}$. Using Remark\,\,\hl{R2.2d}{2.2d}, the local monodromy of the local system $L'_{\xi,\C,\la}$ around $\Yo_{\xi}\cap\Yo_{\xi'_j}$ is given by multiplication by
$$\la^{-\de_{\si,\xi'_j}},$$
\par\nin Note that the multiplicity of the pull-back of $f$ along $\Y_{\xi'_j}$ is $\de_{\si,\1\xi'_j}$, see also Remark\,\,\hl{R2.2e}{2.2e}. By the maximality of $\xi'$, we get the equality of affine subspaces
$$\si\pl\xi'{}^{\perp}=\si\pl\eta_{\si}^{\perp},$$
\par\nin (counting their dimensions). This implies that
$$\de_{\si,\1\xi'}=\de_{\si}.$$
\par\nin The assertion then follows from Remark\,\,\hl{R2.2e}{2.2e} applied to $\xi$ and $\xi'$. Indeed, together with the above equality, it gives that
$$\de_{\si}={\rm GCD}\bl(\de_{\si,\,\xi},\,\de_{\si,\,\xi'_j}\,(j\in[d_{\xi'}{-}d_{\xi}])\br).$$
\par\nin Note that
$$\la^m\ne 1\q\h{if}\q\la\in\mu_{\de_{\si,\xi}}\setminus\mu_{\de_{\si}},\,\,\,m=\de_{\si}.$$
\par\nin Here we can replace the right-hand side of the last condition, that is, $\de_{\si}$ with
$$\de'_{\xi'}:={\rm GCD}\bl(\de_{\si,\,\xi'_j}\,(j\in[d_{\xi'}{-}d_{\xi}])\br),$$
\par\nin using the above equality. Indeed, we have $\la\in\mu_{\de_{\si,\xi}}$ and $\de_{\si}={\rm GCD}(\de_{\si,\,\xi},\de'_{\xi'})$, hence
$$\de_{\si}\1\Z=\de'_{\xi'}\1\Z\,\,\,\,\h{mod}\,\,\,\de_{\si,\,\xi}\1\Z.$$ So the non-triviality of some local monodromy follows from the above calculation of the local monodromies at infinity, since the normal crossing divisors meet at $\Y_{\xi'}$. (Recall that {\it locally\1} the complement of a divisor with normal crossings on a complex manifold has an abelian fundamental group generated freely by loops around local irreducible components, since it is essentially a product of punctured disks.)
We thus get the equality (\hl{2.3.5}{2.3.5}).
\sk
As for the right-hand side of (\hl{2.3.5}{2.3.5}), we can show the following equality in $K_0^{T_s}({\rm MHS})$\,:
\htt{2.3.6}{}
$$\chi_c\bl(\X_{\si},\RR(j_{\si})_*(L_{\si},T_s)\br)=(-1)^{d_{\si}}\bl[H_c^{d_{\si}}\bl(\X_{\si},\RR(j_{\si})_*(L_{\si},T_s)\br)\br].
\leqno(2.3.6)$$
\par\nin Indeed, there is a commutative diagram
$$\begin{array}{cccccccc}\X_{\si}^{\circ}&\buildrel{j_{\si}}\over\longrightarrow&\X_{\si}\\ \,\,\,\downarrow\scriptstyle{\,\jh_{\si}^{\,\circ}}&&\,\,\,\downarrow\scriptstyle{\,\jh_{\si}}\raise5mm\h{}\\ \Xo_{\si}^{\circ}&\buildrel{\jo_{\si}}\over\longrightarrow&\Xo_{\si}\\
\end{array}$$
\par\nin By Remark\,\,\hl{R2.1e}{2.1e}, $\Xo_{\si}^{\circ}$ is an {\it affine\1} variety, and all the inclusions are {\it affine\1} morphisms so that the associated open direct images are $t$-exact functors (see \cite{BBD}). We have moreover the canonical isomorphism
\htt{2.3.7}{}
$$(\jh_{\si})_!\RR(j_{\si})_*L_{\si}=\RR(\jo_{\si})_*(\jh_{\si}^{\,\circ})_!L_{\si},
\leqno(2.3.7)$$
\par\nin using the non-degeneracy condition on $f$. This can be reduced to the case $\Xo_{\si}$ is smooth and $\Xo_{\si}\setminus\X_{\si}$ is a divisor with normal crossing, taking a subdivision of the dual fan of $\si$ (see for instance \cite[8.2]{Da1}) and using the direct image by the morphism from the desingularization. The assertion (\hl{2.3.6}{2.3.6}) then follows from Artin's vanishing theorem \cite{BBD} (see also \cite[2.1.18]{mhp}) together with its dual.
\sk
For $\eta\in\Si'$, let $\Xi'_{\eta}\subset\Xi'$ be the subset consisting of $\xi\in\Xi'$ such that $\eta$ coincides with $\eta_{\xi}$ (see Remark\,\,\hl{R2.1d}{2.1d}), that is, the {\it interior\1} of $\xi$ is contained in that of $\eta$. (In this paper, the interior of a cone means the interior in the smallest affine space containing the cone.) The last condition is equivalent to that $\rho$ induces the morphism $\Y_{\xi}\to\X_{\eta}$. Set
$$I_0:=\{i\in[1,n]\mid e_i\in\eta\},$$
\par\nin where $e_i\in\R^n$ is the $i$\,the unit vector. For $I\subset I_0$ (here $I$ may be $\emptyset$), set
$$\Phi_I:=\msqcup_j\,\Phi_I^j,\q\Phi_I^{\circ}:=\msqcup_j\,\Phi_I^{\circ j},$$
\par\nin with
$$\aligned\Phi_I^j&:=\bl\{\xi\in\Xi'_{\eta}\,\big|\,e_i\in\xi\,\,(\forall\,i\in I),\,d_{\eta}{-}\1 d_{\xi}=j\br\}\,=\,\mcap_{i\in I}\,\Phi_i^j,\\ \Phi_I^{\circ j}&:=\Phi_I^j\setminus\mcup_{i'\in I_0\setminus I}\,\Phi_{i'}^j\,=\,\mcap_{i\in I}\,\Phi_i^j\setminus\mcup_{i'\in I_0\setminus I}\,\Phi_{i'}^j,\raise5mm\h{}\endaligned$$
\par\nin where $\Phi_i^j:=\Phi_{\{i\}}^j$. We have
$$\Xi'_{\eta}=\msqcup_{I\subset I_0}\,\Phi_I^{\circ},$$
\par\nin where $I$ runs over any subsets of $I_0$ including $\emptyset$ (and $\Phi_{I_0}^{\circ}=\Phi_{I_0}$). Put
$$\ep_I:=\msum_j\,(-1)^j|\Phi_I^j|,\q \ep^{\circ}_I:=\msum_j\,(-1)^j|\Phi_I^{\circ j}|.$$
\par\nin Then Theorem\,\,\hl{T1}{1} follows from (\hl{2.3.1}{2.3.1}) and (\hl{2.3.5}{2.3.5}--\hl{2.3.6}{6}) using Proposition\,\,\hl{P2.3}{2.3} and Corollary\,\,\hl{C2.3}{2.3} below.
\par\htt{P2.3}{}\msn
{\bf Proposition\,\,2.3.} {\it In the above notation, we have $\ep_I=1$ for any $I\subset I_0$ including the case $I=I_0$ or $I=\emptyset$.}
\msn
{\it Proof.} Set
$$Y:=\{(y_1,\dots,y_n)\in(\R_{>0})^n\mid\msum_i\,y_i=1\}.$$
\par\nin For $I\subset I_0$, let $Y_I\subset Y$ be the subset defined by $y_i=0$ for $i\notin I$. Let
$Z_I$ be a sufficiently general affine subspace of $Y$ intersecting $Y_I$ transversally. We then get a stratification of an open subset $U_I\subset Z_I$ by
$$Z_{I,\xi}^{\circ}:=\xi^{\circ}\cap Z_I\q(\xi\in\Phi_I),$$
\par\nin (where $\xi^{\circ}$ denotes the interior of $\xi$ in the smallest affine space containing $\xi$) so that
$$U_I=\msqcup_{\xi\,\in\,\Phi_I}\,Z_{I,\xi}^{\circ},\q\h{hence}\q\chi_c(U_I)=\msum_{\xi\,\in\,\Phi_I}\,\chi_c(Z_{I,\xi}^{\circ}).$$
\par\nin Note that $Z_{I,\xi}^{\circ}$ and $U_I$ are contractible. So $\chi_c(Z_{I,\xi}^{\circ})=(-1)^{\dim Z_{I,\xi}^{\circ}}$ by duality \cite{Ve}, and similarly for $U_I$. (Indeed, the dualizing complex on an orientable real manifold is given by the constant sheaf shifted by the real dimension.) The assertion then follows.
\par\htt{C2.3}{}\msn
{\bf Corollary\,\,2.3.} {\it In the above notation, we have $\ep^{\circ}_I=0$ for any $I\subset I_0$ including the case $I=\emptyset$ but not for $I=I_0$.}
\msn
{\it Proof.} We proceed by decreasing induction on $|I|<|I_0|$.
We have the partition
$$\Psi_I=\msqcup_{I\subset J\subset I_0}\,\Psi_J^{\circ},$$
\par\nin where $J$ runs over any subsets of $I_0$ containing $I$ (including the case $J=I$ or $J=I_0$). Hence
$$\ep_I=\msum_{I\subset J\subset I_0}\,\ep_J^{\circ},$$
\par\nin and we have by Proposition\,\,\hl{P2.3}{2.3}
$$\ep_{I_0}^{\circ}=\ep_{I_0}=\ep_I=1.$$
\par\nin By inductive hypothesis, we have $\ep_J^{\circ}=0$ for any $J$ such that $I\,\nsset\,J\,\nsset\,I_0$. (In the case $|I|=|I_0|-1$, this is trivially satisfied, since there is no such $J$.) We thus get $\ep_I^{\circ}=0$. So the assertion follows by induction.
\par\htt{R2.3}{}\msn
{\bf Remark\,\,2.3.} Related to $\psi_{\si}$ in the case $\si=\emptyset$ in (\hl{2}{2}), we have in the $\si$ {\it simplicial\1} case
$$H^k\bl(\Xo_{\si}\setminus{}\,\Xo_{\si}^{\rm pr},\Q\br)=0\q(k\ne 0,d_{\si}),$$
\par\nin with $\Xo_{\si}^{\rm pr}:=\Xo_{\si}\cap\X_0^{\rm pr}$. Here $H^0(\Xo_{\si}\setminus\,\Xo_{\si}^{\rm pr},\Q)=\Q$, and $H^{d_{\si}}(\Xo_{\si}\setminus{}\,\X_{\si}^{\rm pr},\Q)$ is pure of weight $d_{\si}{+}1$, using the Gysin sequence
$$\to H^{k-2}(\Xo_{\si}^{\rm pr},\Q)(-1)\to H^k(\Xo_{\si},\Q)\to H^k(\Xo_{\si}\setminus{}\,\Xo_{\si}^{\rm pr},\Q)\to
$$
\par\nin and the weak Lefschetz property. Note that $\Xo_{\si}$ is a $V$-manifold which is homeomorphic to $\PP^{d_{\si}}$ (since $\si$ is simplicial), and $\Xo_{\si}^{\it pr}$ is also a $V$-manifold (since $f$ is non-degenerate).
Using a spectral sequence, we then get that
$$\Gr_0^WH_c^{d_{\si}}(\X_{\si}^{\circ},\Q)=\Q\q(\si\ne\emptyset).$$
\par\nin This shows that $\psi_{\emptyset}$ is necessary in (\hl{2}{2}). Indeed, $\Gr_0^WH^{n-1}(\Ff,\Q)_1$ does not vanish without it, since $\msum_{i\ges 0}\,(-1)^i\1s_i=1$ with $s_i$ the number of $i$-dimensional faces of $\Gf$.
\par\htt{2.4}{}\msn
{\bf 2.4.~Proofs of Corollary\,\,\hl{C1}{1} and Theorem\,\,\hl{T2}{2}.} We may assume that $f$ is {\it typical\1} by Remarks\,\,\hl{R2.4b}{2.4b}--\hl{R2.4c}{c} below. (Here typical means that ${\rm Supp}_{(x)}f$ consists of the vertices of $\Gp(f)$.) Then the assertions follow from Theorem\,\,\hl{T1}{1} using Lemma\,\,\hl{L2.4}{2.4} and Proposition\,\,\hl{P2.4}{2.4} below.
\par\htt{L2.4}{}\msn
{\bf Lemma\,\,2.4.} {\it Assume $\si\les\Gf$ is a simplex. Set $\Xo^{\circ}_{\si}:=\Xo_{\si}\setminus\X_0^{\rm pr}$ with $\jh^{\,\circ}_{\si}:\X_{\si}^{\circ}\into\Xo^{\circ}_{\si}$ the inclusion, and
$$\Lo_{\si}:=(\jh^{\,\circ}_{\si})_*L_{\si}.$$
\par\nin There are isomorphisms constructible sheaves on $\Xo^{\circ}_{\si}$ for $\emptyset\ne\tau<\si:$
\htt{2.4.1}{}
$$\Lo_{\si}|_{\Xo^{\circ}_{\tau}}=\Lo_{\tau}.
\leqno(2.4.1)$$
\par\nin Moreover $\Lo_{\si}$ is identified with the {\it intermediate direct image\1} of $L_{\si}$ up to a shift of complex.}
\msn
{\it Proof.} This follows from an argument similar to \hl{2.3}{2.3}. Indeed, the assertion is reduced to the normal crossing case by using the morphism $\rho$ and choosing $\xi_1,\xi_2\in\Xi'$ such that
$$\xi_1\subset\xi_2,\,\,\,\xi_1\subset\eta_{\si},\,\,\,\xi_2\subset\eta_{\tau},\,\,\,d_{\xi_1}=d_{\eta_{\si}},\,\,\,d_{\xi_2}=d_{\eta_{\tau}}.$$
\par\nin The assertion then follows by using Remarks\,\,\hl{R2.2d}{2.2d}--\hl{R2.2e}{e}. (This is closely related to Remark\,\,\hl{R2.1g}{2.1g}.)
Note that the assumption on $\si$ implies that the $\Xo_{\tau}$ are $V$-manifolds so that a constant sheaf is an intersection complex up to shift for $\tau\les\si$, see also the proof of Proposition\,\,\hl{P2.4}{2.4} and Remark\,\,\hl{R2.4a}{2.4a} below. This finishes the proof of Lemma\,\,\hl{L2.4}{2.4}.
\par\htt{P2.4}{}\msn\vbox{\nin
{\bf Proposition\,\,2.4.} {\it Assume $\si\les\Gf$ is a simplex, and moreover $f_{\si}$ is typical. Set
$$\Lo_{\si,\C,\la}:={\rm Ker}(T_s\mi\la)\,\subset\,\C\otimes_{\Q}\Lo_{\si}.$$
\par\nin We have the isomorphisms for $p=[\al]$, $\la=\exp(2\pi i\al)$ with $\al\in\Q$\,$:$
\htt{2.4.2}{}
$$\Gr_F^pH^{d_{\si}}\bl(\Xo_{\si},\RR(\jo_{\si})_*\Lo_{\si,\C,\la}\br)=B_{\si,\1 d_{\si}+1-\al},
\leqno(2.4.2)$$
\par\nin where $\jo_{\si}:\Xo_{\si}^{\circ}\into\Xo_{\si}$ is a natural inclusion, and $B_{\si}$ is as in the introduction $($defined before Corollary\,\,$\hl{C1}{1})$.}}
\msn
{\it Proof.} Let $v_1,\dots,v_r$ be the vertices of $\si$. Let $\de_{\si}$, $V(\si)$ be as in (\hl{1}{1}). Set $E':=V(\si)\cap\Z^n$. Let $E\subset V(\si)$ be the free abelian subgroup generated by $v'_k:=\tfrac{\!\!1}{\de_{\si}}v_k$ ($k\in[1,r]$), where $r=d_{\si}+1$. Then $E\supset E'$, and Remark\,\,\hl{R2.4a}{2.4a} below applies. Let $x_i\in\C[E]$ corresponding to $v'_k\in E$. Then we have
$$\C[x_1,\dots,x_r]^G=A_{\si},$$
\par\nin since each monomial is stable by the action of $G$ (up to constant multiple). Moreover, there is $h=\msum_{i=1}^r\,a_i\1x_i^d$ ($a_i\ne 0$) in the left-hand side which is identified with $f_{\si}$ in the right-hand side. Set
$$B^{\circ}_{k/d}:=\{\nu\in\Z^r\mid 0<\nu_i<d\,(\forall\, i),\,\,|\nu|=k\}.$$
\par\nin It is well known (as a consequence of \cite{Br}, \cite{SS}, \cite{Va2}) that we have a canonical isomorphism
$$\iota_h:B^{\circ}_{k/d}\simto\Gr_F^pH^{r-1}(F_{\!h},\C)_{\la},$$
\par\nin for $p=[n\mi k/d]$, $\la=\exp(-2\pi ik/d)$, where $F_{\!h}:=h^{-1}(1)\subset\C^r$. This isomorphism can be defined by
$$\iota_h(g):=\bl[\bl((g\1\ddd x/x)/\ddd h\br)|_{F_{\!h}}\br]\in H^{r-1}_{\rm DR}(F_{\!h},\C)\q(g\in B^{\circ}_{k/d}),$$
\par\nin using the analytic de Rham cohomology of the stein manifold $F_{\!h}$, where $\ddd x/x:=\bigwedge_{i=1}^r\ddd x_i/x_i$. (Note that $g/x\in\C[x_1,\dots,x_r]$.) The isomorphism $\iota_h$ is compatible with the action of $G$, since $\ddd x/x$ and $h$ are $G$-invariant. So the assertion follows by restricting the isomorphism to the $G$-invariant part, since the local system $L_{\si}$ and the constructible sheaf $\Lo_{\si}$ can be obtained also by using the direct images of the constant sheaves on the second and first terms of (\hl{2.1.6}{2.1.6}). (Note that (\hl{2.4.1}{2.4.1}) corresponds to that $f_{\tau}$ is a ``restriction" of $f_{\si}$ for $\tau<\si$.) The index $d_{\si}\pl 1\mi\al$ on the right-hand side of (\hl{2.4.2}{2.4.2}) comes from the difference of the two definitions of spectral numbers as is explained after (\hl{1.1.1}{1.1.1}) (using the symmetry). This finishes the proof of Proposition\,\,\hl{P2.4}{2.4}.
\par\htt{R2.4a}{}\msn
{\bf Remark\,\,2.4a.} Let $E$ be a free abelian group with $E'\subset E$ a subgroup such that $E/E'$ is finite. There is a decomposition
$$E/E'\cong\mprod_{j=1}^b\,\mu_{a_j}\q(a_j\ges 2),$$
\par\nin with $\mu_{a_j}$ as in (\hl{2.2.3}{2.2.3}).
Let $\phi_j:E\to\mu_{a_j}$ be the composition of the projection $E\onto E/E'$ and the projection to the $j$\1th factor.
\sk
Consider the group ring $\C[E]$. We denote by $x^e$ the element of $\C[E]$ corresponding to $e\in E$ so that $x^{e+e'}=x^ex^{e'}$.
(If we choose free generators $e_1,\dots,e_r$ of $E$, then $\C[E]$ is identified with $\C[x_1,\dots,x_r]\bl[\tfrac{1}{x_1\cdots x_r}\br]$ where $x_k:=x^{e_k}$.) We have an action of a generator $\rho_j$ of $\mu_{a_j}$ on $\C[E]$ such that
$$\rho_j(x^e)=\phi_j(e)\1x^e\q(e\in E).$$
\par\nin This is extended to the action of $G:=\mprod_{j=1}^b\,\mu_{a_j}$ on $\C[E]$ so that
$$\C[E]^G=\C[E'].$$
\par\nin \par\htt{R2.4b}{}\msn
{\bf Remark\,\,2.4b.} If there is a deformation $\{f_u\}_{u\in\De^r}$ such that $\Gp(f_u)$ is independent of $u$ and each $f_u$ is Newton non-degenerate, then we have a simultaneous embedded resolution of the $f_u$ by taking a smooth subdivision of the dual fan of $\Gp(f_u)$, see for instance \cite[8.2]{Da1}. This implies that the spectrum of $f_u$ is independent of $u$ (similarly for spectral pairs in the isolated singularity case).
\par\htt{R2.4c}{}\msn
{\bf Remark\,\,2.4c.} For a finite number of lattice points $v^{(j)}\in\Z^n$ ($j\in[1,r]$), there is a non-empty Zariski-open subset $U_v\subset\C\1^r$ such that, for $(a_j)\in U_v$, we have the smoothness of the hypersurface
$$Z_{v,a}:=\bl\{\1(x_i)\in(\C^*)^n\,\big|\,\msum_{j=1}^r\,a_j\1x^{v^{(j)}}=0\1\br\}.$$
\par\nin Indeed, the hypersurface defined by the above equation in $(\C^*)^n{\times}U_v$ is smooth (using the derivation by $a_i$). So the assertion follows from a well-known Bertini type theorem.
\par\htt{R2.4d}{}\msn
{\bf Remark\,\,2.4d.} The compatibility of Corollary\,\,\hl{C1}{1} with \cite{Va1} can be verified by using an argument similar to an old version of \cite[1.5]{JKSY2}.
\par\htt{R2.4e}{}\msn
{\bf Remark\,\,2.4e.} We can show that Corollary\,\,\hl{C1}{1} implies \cite[Theorem 2]{JKSY2}, which claims that for a non-degenerate function $f$ of 3 variables with $\Gp(f)$ simplicial and intersecting any coordinate plane of $\R^3$, we have
\htt{2.4.3}{}
$$\aligned\Sp'_f(t)&=\msum_{\si\in\CFfin}\,\bl(\msum_{j=0}^{2-d_{\si}}\,t^j\br)\,q_{\si}(t)+\bl|\CFfin^0\br|\,\bl(t{+}t^2\br)\\&\q{}+{}\msum_{\si\in\CFf^0}\,\bl(\gat_{\si}{-}3\br)\1q_{\si}(t)\1t\\&\q{}-\msum_{\si\in\CFf}\,m_{\si}\bl(\msum_{j=0}^{1-d_{\si}}\,t^j\br)\,q_{\si}(t)-\msum_{\si\in\CFf^0}\,m_{\si}\1t.\endaligned
\leqno(2.4.3)$$
\par\nin Here $m_{\si}$ for $\si<\Gp(f)$ is the number of $(d_{\si}{+}1)$-dimensional faces of $\Gp(f{+}\ell\1^r)$ ($r\,{\gg}\,0)$ which is the convex hull of $\si\cup\{r\1\ee_i\}$ for some $i\in[1,3]$ and is not contained in any coordinate plane with $\ee_i$ the $i$\1th unit vector and $\ell$ a sufficiently general linear function. We denote by $\gat_{\si}$ the number of $2$-dimensional $($not necessarily \h{\it compact\1$)$} faces of $\Gp(f{+}\ell\1^r)$ ($r\gg 0)$ containing $\si\in\CFf^0$.
\sk
Indeed, let $\CFfsb^j$, $\CFfdb^j$ be the subset of $\CFfb^j$ consisting of $\tau$ such that $\si^{(i)}\cup\tau$ is a face of $\Gp(f{+}\ell\1^r)$ ($r\gg 0$) for exactly one and two $i\in[1,3]$ respectively. Put
$$\aligned\CFfnb^j&:=\CFf^j\setminus\CFfb^j,\\\CFf^{j,k}&:=\bl\{\si\in\CFf^j\,\big|\,k(\si)=k{+}1\br\}.\endaligned$$
\par\nin Note that $\CFf^{j,2}=\CFfin^j$. Set
$$\CFfb^{j,k}:=\CFfb^j\cap\CFf^{j,k},$$
\par\nin and similarly for $\CFfsb^{j,k}$, etc. Put
$$n_{\si,k,d}:=\#\{\1\tau\les\Gf\mid\tau\ges\si,\,k(\tau)=k{+}1,\,d(\tau)=d{+}1\1\}.$$
\par\nin Note that
$$\msum_k\,n_{\si,k,1}+m'_{\si}=\gat_{\si}\q\h{if}\q\si\in\CFf^0.$$
\par\nin Here $m'_{\si}$ is defined like $m_{\si}$ without assuming the last condition (about the non-inclusion in coordinate planes).
\sk
We can calculate these $n_{\si,k,d}$ and the combinatorial polynomials $r_{\si}(t)$ (using pictures as in Remark\,\,\hl{R2.4f}{2.4f} below if necessary) as follows.
\msn
Case 0a: $\,\si\in\CFfnb^{0,2}$ ($m_{\si}=0$).
$$\aligned&n_{\si,2,0}=1,\q n_{\si,2,1}=\gat_{\si},\q n_{\si,2,2}=\gat_{\si},\q (d_{\si}=0),\\&r_{\si}(t)=
(t{-}1)^2+\gat_{\si}(t{-}1)+\gat_{\si}=(\gat_{\si}{-}3)t+(1{+}t{+}t^2).\endaligned$$
\par\nin \skn
Case 0b: $\,\si\in\CFfsb^{0,2}$ ($m_{\si}=1$).
$$\aligned&n_{\si,2,0}=1,\q n_{\si,2,1}=\gat_{\si}-1,\q n_{\si,2,2}=\gat_{\si}-2,\\&r_{\si}(t)=
(t{-}1)^2+(\gat_{\si}{-}1)(t{-}1)+(\gat_{\si}{-}2)\\&\q\q\,=(\gat_{\si}{-}3)t+(1{+}t{+}t^2)-(1{+}t).\endaligned$$
\par\nin \skn
Case 0b$'$: $\,\si\in\CFfdb^{0,2}0$ ($m_{\si}=2$).
$$\aligned&n_{\si,2,0}=1,\q n_{\si,2,1}=\gat_{\si}-2,\q n_{\si,2,2}=\gat_{\si}-4,\\&r_{\si}(t)=
(t{-}1)^2+(\gat_{\si}{-}2)(t{-}1)+(\gat_{\si}{-}4)\\&\q\q\,=(\gat_{\si}{-}3)t+(1{+}t{+}t^2)-2(1{+}t).\endaligned$$
\par\nin \skn
Case 0c: $\,\si\in\CFfnb^{0,1}$ ($m_{\si}=0$).
$$\aligned&n_{\si,1,0}=1,\q n_{\si,2,1}=\gat_{\si}-2,\q n_{\si,1,1}=2,\q n_{\si,2,2}=\gat_{\si}-1,\\&r_{\si}(t)=
-(t{-}1)+(\gat_{\si}{-}2)(t{-}1)-2+(\gat_{\si}{-}1)=(\gat_{\si}{-}3)t.\endaligned$$
\par\nin \skn
Case 0d: $\,\si\in\CFfsb^{0,1}$ with $m_{\si}=1$.
$$\aligned&n_{\si,1,0}=1,\q n_{\si,2,1}=\gat_{\si}-3,\q n_{\si,1,1}=2,\q n_{\si,2,2}=\gat_{\si}-3,\\&r_{\si}(t)=
-(t{-}1)+(\gat_{\si}{-}3)(t{-}1)-2+(\gat_{\si}{-}3)=(\gat_{\si}{-}3)t-(t{+}1).\endaligned$$
\par\nin \skn
Case 0d$'$: $\,\si\in\CFfsb^{0,1}$ with $m_{\si}=0$.
$$\aligned&n_{\si,1,0}=1,\q n_{\si,2,1}=\gat_{\si}-2,\q n_{\si,1,1}=1,\q n_{\si,2,2}=\gat_{\si}-2,\\&r_{\si}(t)=
-(t{-}1)+(\gat_{\si}{-}2)(t{-}1)-1+(\gat_{\si}{-}2)=(\gat_{\si}{-}3)t.\endaligned$$
\par\nin \skn
Case 0d$''$: $\,\si\in\CFfdb^{0,1}$ ($m_{\si}=0$).
$$\aligned&n_{\si,1,0}=1,\q n_{\si,2,1}=\gat_{\si}-2,\q n_{\si,1,1}=0,\q n_{\si,2,2}=\gat_{\si}-3,\\&r_{\si}(t)=
-(t{-}1)+(\gat_{\si}{-}2)(t{-}1)-0+(\gat_{\si}{-}3)=(\gat_{\si}{-}3)t.\endaligned$$
\par\nin \skn
Case 0e: $\,\si\in\CFf^{0,0}$ ($m_{\si}=0$).
$$\aligned&n_{\si,0,0}=1,\q n_{\si,2,1}=\gat_{\si}-3,\q n_{\si,1,1}=2,\q n_{\si,2,2}=\gat_{\si}-2,\\&r_{\si}(t)=
1+(\gat_{\si}{-}3)(t{-}1)-2+(\gat_{\si}{-}2)=(\gat_{\si}{-}3)t.\endaligned$$
\par\nin \skn
Case 1a: $\,\si\in\CFfnb^{1,2}$ ($m_{\si}=0$).
$$n_{\si,2,1}=1,\q n_{\si,2,2}=2,\q r_{\si}(t)=(t{-}1)+2=t{+}1.$$
\par\nin \skn
Case 1b: $\,\si\in\CFfb^{1,2}$ ($m_{\si}=1$).
$$n_{\si,2,1}=1,\q n_{\si,2,2}=1,\q r_{\si}(t)=(t{-}1)+1=(t{+}1)-1.$$
\par\nin \skn
Case 1c: $\,\si\in\CFfnb^{1,1}$ ($m_{\si}=0$).
$$n_{\si,1,1}=1,\q n_{\si,2,2}=1,\q r_{\si}(t)=-1+1=0.$$
\par\nin \skn
Case 1d: $\,\si\in\CFfb^{1,1}$ ($m_{\si}=1$).
$$n_{\si,1,1}=1,\q n_{\si,2,2}=0,\q r_{\si}(t)=-1.$$
\par\nin \ms
We now get the equality between the partial sum of (\hl{3}{3}) over the non-empty faces $\si\les\Gf$ and the right-hand side of (\hl{2.4.3}{2.4.3}) with second and last terms deleted. The coincidence between $r_{\emptyset}(t)$ and the sum of the second and last terms of the right-hand side of (\hl{2.4.3}{2.4.3}) is reduced to Remark\,\,\hl{R2.4f}{2.4f} just below.
\par\htt{R2.4f}{}\msn
{\bf Remark\,\,2.4f.} For a triangulation of a triangle $\Gamma$, let $a_0,a_1$ be the number of vertices and edges of the triangulation which are not contained in the boundary of $\Gamma$. Let $b_0$ be the number of vertices lying on {\it smooth\1} points of $\dd\1\Gamma$. Then we have
$$3\1a_0+b_0=a_1.$$
\par\nin \sk
For instance, if the triangulation is given by
$$\setlength{\unitlength}{.4cm}
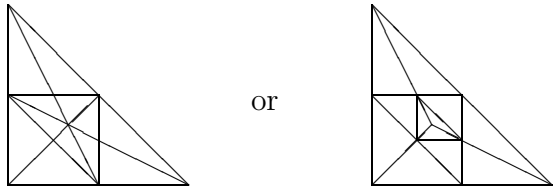
\begin{picture}(6,6)
\put(0,0){\line(0,1){6}}
\put(0,0){\line(1,0){6}}
\put(6,0){\line(-1,1){6}}
\put(6,0){\line(-2,1){6}}
\put(0,6){\line(1,-2){3}}
\put(0,0){\line(1,1){3}}
\put(0,3){\line(1,-1){3}}
\put(0,3){\line(1,0){3}}
\put(3,0){\line(0,1){3}}
\end{picture}
\q\q\raise1cm\h{or}\q\q\q
\begin{picture}(6,6)
\put(0,0){\line(0,1){6}}
\put(0,0){\line(1,0){6}}
\put(6,0){\line(-1,1){6}}
\put(6,0){\line(-2,1){4}}
\put(0,6){\line(1,-2){2}}
\put(0,0){\line(1,1){2}}
\put(0,3){\line(1,-1){3}}
\put(0,3){\line(1,0){3}}
\put(3,0){\line(0,1){3}}
\put(1.5,1.5){\line(0,1){1.5}}
\put(1.5,1.5){\line(1,0){1.5}}
\put(1.5,3){\line(1,-1){1.5}}
\end{picture}
$$
\par\nin we have $(a_0,b_0,a_1)=(4,3,15)$. Here a triangulation means that if the intersection of two triangles $\ga_1,\ga_2$ is non-empty, then it coincides with one side of $\ga_i$ for $i=1,2$.
\sk
The above statement can be verified by induction on $a_0$ dividing $\Ga$ into two parts. Here we may assume $a_0>0$, since the case $a_0=0$ is easily shown by induction on $a_1$ or $b_0$ (removing an appropriate edge).
We choose a general line $\ell$ passing through a vertex $A$ of $\Ga$ and a point $B$ on the boundary of $\Ga$ which is not a vertex of $\Ga$. We assume $\ell$ does not contain any vertex of triangulation (except for $A$). Let $\Ga_1,\Ga_2$ be the triangles such that $\Ga_1\cup\Ga_2=\Ga$, $\Ga_1\cap\Ga_2=\Ga\cap\ell$. We may assume that $\Ga_i$ contains at least one vertex of triangulation for $i=1,2$ (replacing $A$ if necessary). Let $\ga_j$ ($j\in[0,r]$) be the triangles intersecting $\ell$. Changing the order of the $\ga_j$ if necessary, we may assume that there are intersection points $P_j$ ($j\in[1,r]$) of $\ell$ with edges of triangulation such that $P_j,P_{j+1}\in\ga_j$ ($j\in[0,r]$), where $P_0=A$, $P_{r+1}=B$. We then get triangulations of $\Ga_1,\Ga_2$ by adding an edge joining $P_j$ and a vertex of $\ga_j$ for each $j\in[1,r]$. Define $a^{(i)}_0,b^{(i)}_0,a^{(i)}_1$ for $\Ga_i$ as above ($i=1,2$). Since the edges intersecting $\ell$ are divided into two parts, we see that
$$a^{(1)}_0\pl a^{(2)}_0=a_0,\q b^{(1)}_0\pl b^{(2)}_0=b_0\pl 2\1 r,\q a^{(1)}_1\pl a^{(2)}_1=a_1\pl 2\1r.$$
\par\nin The assertion is thus reduced to the case $a_0=1$ by induction. We may then assume that the unique vertex in the interior of $\Ga$ is contained in all the edges by deleting (or replacing) inductively certain edges which are a side of a triangle containing a vertex of $\Ga$. The claim is further reduced to the case $b_0=0$ by removing edges inductively. We then get $a_1=3$. So the assertion follows.
\sk
For the application to the problem in Remark\,\,\hl{R2.4e}{2.4e} above, all the triangles containing certain vertices of $\Ga$ must be deleted (since $f$ has non-isolated singularities). The  multiplicity of each vertex on the new boundary of the remaining area then becomes $3\mi m_{\si}$ if $\si\in\CFf^0$ corresponds to the vertex (since the edge joining this vertex with a vertex of $\Ga$ is removed). Here the triangulation by obtained by taking the intersection of the cones of $\si\les\Gf$ with the hypersurface $\msum_i\,\nu_i=1$. (This triangulation cannot recover the toric variety associated with $\Si$, since crucial information is lost by passing to the cones of $\si$.)
\par\htt{2.5}{}\msn
{\bf 2.5.~Proof of Corollary\,\,\hl{C3}{3}.} By Remarks\,\,\hl{R2.4b}{2.4b}--\hl{R2.4c}{c} just above, we may assume $f$ is typical as in \hl{2.4}{2.4}. Using the intermediate direct image explained after (\hl{2.2.7}{2.2.7}), the assertion on $n_{\la,k}$ for $\la\ne 1$, $k=n$ or $n{-}1$ follows from the {\it estimate of weights\1} as in \cite[4.5.2]{mhm} and the monodromical property of the weight filtration as in Remark\,\,\hl{R1.2b}{1.2b} together with Proposition\,\,\hl{P2.4}{2.4} (since any 1-dimensional polytope is a simplex). Indeed, the estimate of weights implies that ${\rm wt}(\psi_{\si,\1\ne1})\les d_{\si}$, and hence
\htt{2.5.1}{}
$${\rm wt}\bl((1-\theta)^{k(\si)-d(\si)}\1\psi_{\si,\1\ne1}\br)\les 2n\mi2\mi d_{\si},
\leqno(2.5.1)$$
\par\nin since $k(\si)-d(\si)\les n-1-d_{\si}$ (and ${\rm wt}\,\theta=2$). Here the decomposition $\psi_{\si}=\psi_{\si,1}+\psi_{\si,\1\ne1}$ (induced by (\hl{2.2.4}{2.2.4})) is used, and wt means weights.
\sk
The assertion on $n_{1,n-1}$ also follows from Theorem\,\,\hl{T1}{1} and Proposition\,\,\hl{P2.4}{2.4} using the estimation
\htt{2.5.2}{}
$${\rm wt}\bl((1-\theta)^{k(\si)-d(\si)}\1\psi_{\si,1}\br)\les 2n\mi1\mi d_{\si},
\leqno(2.5.2)$$
\par\nin (which holds since the closure of $\X_{\si}^{\rm pr}$ in an appropriate compactification of $\X_{\si}$ is a smooth hypersurface).
This finishes the proof of Corollary\,\,\hl{C3}{3}.
\par\htt{R2.5a}{}\msn
{\bf Remark\,\,2.5a.} For the proof of the assertion on $n_{1,n-1}$, we can also use the equality
\htt{2.5.3}{}
$$\dim\Gr^W_2H^1(\PP\setminus Z)=|Z|-1,
\leqno(2.5.3)$$
\par\nin for a non-empty finite subset $Z\subset\PP^1$ instead of Proposition\,\,\hl{P2.4}{2.4}.
\par\htt{R2.5b}{}\msn
{\bf Remark\,\,2.5b.} The assertion for $\la\ne 1$, $k=n{-}1$ in Corollary\,\,\hl{C3}{3} is equivalent to the equality
\htt{2.5.4}{}
$$n_{\la,n-1}=\msum_{\si\les\Gf,\,d_{\si}=1}\,(q_{\si,\be}+q_{\si,2-\be}),
\leqno(2.5.4)$$
\par\nin for $\be\in(0,1)$ with $\la=\exp(\pm 2\pi i\be)$ (here $q_{\si}(t)=\msum_{\al\in\Q}\,q_{\si,\al}\1t^{\al}$), see also \cite{Sta} and the references quoted there as relevant results. We have $n_{\la,k}=n_{\overline{\la},k}$, since the monodromy is defined over $\Q$.
\par\htt{2.6}{}\msn
{\bf 2.6.~Proof of Corollary\,\,\hl{C2}{2}.} Let $-\al_f$ be the maximal root of the Bernstein-Sato polynomial $b_{f,0}(s)$. By \cite[Corollary 3.3]{mic} (using implicitly \cite[(1.3.4)]{hi}), we have
\htt{2.6.1}{}
$$\al_f\ges c^{-1}.
\leqno(2.6.1)$$
\par\nin So it is sufficient to show that $b_{f,0}(-c^{-1})=0$. Using the projections $\R^n\to\R^{n-1}$ together with Remark\,\,\hl{R2.6a}{2.6a} below, the assertion is reduced to the case $(c,\dots,c)\in\Gf$.
\sk
It is then further reduced to the assertion that $c^{-1}$ is a spectral number of the Steenbrink spectrum $\Sp_f(t)=\Sp'_f(t^{-1})\1t^n$. Indeed, by \cite[Theorem 2]{bfun}, $-\al$ is a root of $b_{f,0}(s)$ for $\al\in\Q_{>0}$ if
\htt{2.6.2}{}
$$\Gr_{\Pt}^pH^{n-1}(\Ff,\C)_{\la}\ne0\q\q(\la=e^{-2\pi i\al},\,\,p=[n{-}\al]).
\leqno(2.6.2)$$
\par\nin On the other hand, the spectral numbers $\al$ of $\Sp_f(t)$ contained in $(0,1)$ are given by (\hl{2.6.2}{2.6.2}) with $\Pt$ replaced by $F$, since we have $p=[n\mi\al]=n{-}1$, and
\htt{2.6.3}{}
$$F^{n-1}H^k(\Ff,\C)=0\q\h{if}\q k<n{-}1,
\leqno(2.6.3)$$
\par\nin see Remark\,\,\hl{R2.6a}{2.6a} below. Moreover the saturated pole order filtration $\Pt$ contains the Hodge filtration $F$ on $H^{n-1}(\Ff,\C)$, and $\Pt^n=0$, see {\it loc.\,cit.} So the above reduction is proved.
\sk
By Corollary\,\,\hl{C1}{1}, the {\it maximal\1} spectral number of $\Sp'_f(t)$ is given by $n\mi c^{-1}$, using the symmetry of the $q_{\si}(t)$. Indeed, $f$ is simplicial, the leading coefficient of $r_{\si}(t)$ is 1 for the minimal dimensional $\si\les\Gf$ with $(c,\dots,c)\in\si$, and $c^{-1}<1$. So the minimal spectral number of $\Sp_f(t)$ is $c^{-1}$. This finishes the proof of Corollary\,\,\hl{C2}{2}.
\par\htt{R2.6a}{}\msn
{\bf Remark\,\,2.6a.} If $f$ is Newton non-degenerate, then so is the restriction $f_{(i)}$ of $f$ to $x_i=a_i$ with $|a_i|$ sufficiently small in the case $f_{(i)}$  has a singularity. This can be verified by using a smooth subdivision of the dual fan. Indeed, it implies the {\it stratified smoothness\1} of the pull-back of $f_{(i)}$ when $|a_i|$ is very small, see also \cite[Remark 1.7]{JKSY1}.
\par\htt{R2.6b}{}\msn
{\bf Remark\,\,2.6b.} The mixed Hodge structure on the vanishing cohomology $H^{n-1+k}(\Ff,\Q)$ is defined by applying the cohomological functor $H^ki_0^*$ to the nearby mixed Hodge module of $f$. Here $i_0:\{0\}\into X$ is the inclusion, and $H^ki_0^*$ can be defined by iterating the mapping cones of functors
\htt{2.6.4}{}
$$C({\rm can}_j:\psi_{x_j,1}\to\varphi_{x_j,1})\q(j\in[1,n]),
\leqno(2.6.4)$$
\par\nin with $x_j$ local coordinates, see \cite{mhm}. The morphism ${\rm can}_j$ in (\hl{2.6.4}{2.6.4}) is induced by $\dd_{x_j}$, and the Hodge filtration $F$ on $\psi$ is shifted by 1 so that $F$ is preserved by ${\rm can}_j$, see \cite{mhp}. This shift implies that
\htt{2.6.5}{}
$$F^{n-k+1}H^{n-k}(\Ff,\C)=0\q(k\ges1).
\leqno(2.6.5)$$
\par\nin \par\htt{R2.6c}{}\msn
{\bf Remark\,\,2.6c.} It does not seem easy to generalize the argument in the proof of Corollary\,\,\hl{C2}{2} to the case $c<1$ where the situation is quite different; consider for instance the case $f=x^a\pl yz$ (but not $x^a\pl y^2\pl z^2$) for $a>2$, where $\Sp_f(t)=r_{\si}(t)\1q_{\si}(t)$ for $\si=\{(a,0,0)\}$ with $r_{\si}(t)=t$, $q_{\si}(t)=t^{1/a}+\cdots+t^{(a-1)/a}$.
\bsn\htt{A}{}\par
\vbox{\centerline{\bf Appendix. Descent theorem for motivic nearby fibers}
\bsn
In this Appendix we show the descent theorem for motivic nearby fibers.}
\par\htt{A.1}{}\msn
{\bf A.1.~Relative Grothendieck ring.} In this paper a variety means a separated scheme of finite type over $\C$. We assume further it is {\it reduced,} but not necessarily irreducible (since we consider Grothendieck rings of varieties). Moreover we consider only {\it closed points\1} of varieties. So it is an algebraic variety over $\C$ in the sense of Serre \cite[Section 34]{Se}.
The analytic spaces $\X_0^{\rm ex}$, etc.\ constructed in Section 2 will be viewed as varieties (although the same notation is used), since they are naturally defined algebraically.
\sk
For a variety $S$, let $K_0^{\muh}({\rm Var}_S)[\Lb^{-1}]$ be the relative Grothendieck ring of varieties endowed with a good $\muh$-action over $S$ with $\Lb$ inverted, see \cite[2.4]{DL}, \cite[3.4]{Loe}, \cite[2.2]{Bi}, etc. Here the $\muh$-action on $S$ is trivial, and good action means that every orbit is contained in an affine open subset. The multiplicative structure is defined by fiber product over $S$.
(Note that the fiber product of affine open subvarieties of two varieties over $S$ is the intersection of their fiber product over $\C$ with the inverse image of the diagonal of $S{\times}S$, and closed subvarieties of affine varieties are affine.)
Recall that $\Lb$ is the class of the trivial affine bundle $[\Ab^1{\times}S]$ endowed with the trivial $\muh$-action (here $\Ab^1$ denotes 1-dimensional affine space). Note that $\muh$ is the projective limit of $\mu_m$ ($m\in\Z_{>0}$), and $K_0^{\muh}({\rm Var}_S)$ is the inductive limit of $K_0^{\mu_m}({\rm Var}_S)$, so every object has an $\mu_m$-action for some sufficiently divisible positive integer $m$.
\sk
By definition $K_0^{\muh}({\rm Var}_S)$ is a quotient group of the free group generated by the $[Y]$ with $Y$ isomorphism classes of varieties with good $\muh$ action over $S$. This free group is divided by the subgroup generated by
\htt{A.1.1}{}
$$[Y\setminus Z]+[Z]-[Y],
\leqno{\rm(A.1.1)}$$
\par\nin for any variety $Y$ with a good $\muh$-action over $S$ and any closed subvariety $Z\subset Y$ stable by the $\muh$-action. Note that this subgroup is an ideal (since closed immersions are stable by base change).
Hence $K_0^{\muh}({\rm Var}_S)$ is a ring. Its localization by $\Lb$ is denoted by $K_0^{\muh}({\rm Var}_S)[\Lb^{-1}]$.
\sk
Let $K_0^{\muh}({\rm Var}_S)^{\sim}$ be the quotient ring of $K_0^{\muh}({\rm Var}_S)[\Lb^{-1}]$ by the subgroup generated by
\htt{A.1.2}{}
$$[A]-[A'],
\leqno{\rm(A.1.2)}$$
\par\nin with $A\to Y$, $A'\to Y$ affine bundles of the same rank over a same variety $Y$ over $S$, which are endowed with $\mu_m$-actions compatible with a (same) $\mu_m$-action on $Y$ over $S$ for some $m\in\Z_{>0}$, see \cite[3.4]{Loe}. This subgroup is an ideal (since affine bundles are stable by base change).
\par\htt{RA.1a}{}\msn
{\bf Remark\,\,A.1a.} In (\hl{A.1.2}{A.1.2}) we may assume that $A$, $A'$ are {\it products\1} of affine spaces with a same variety $Y$ over $S$ by using a sufficiently fine partition of $Y$ (since we consider them in the Grothendieck ring), although we cannot assume that the product structures are compatible with the $\mu_m$-actions, that is, the actions are {\it diagonal\1} on the products.
\par\htt{RA.1b}{}\msn
{\bf Remark\,\,A.1b.} In the case of varieties over $\C$, the $\muh$-action can be replaced with an action of $T_s$ of finite order on varieties over $S$ by choosing $\sqrt{-1}\in\C^*$. Indeed, this choice gives an isomorphism
\htt{A.1.3}{}
$$\Z/\Z m\ni a\,\mapsto\,\exp(2\pi\sqrt{-1}\1a/m)\in\mu_m,
\leqno{\rm(A.1.3)}$$
\par\nin such that the {\it canonical generator} $\,\exp(2\pi\sqrt{-1}/m)\in\mu_m$ (with $a=1$) is sent to the canonical generator $\exp(2\pi\sqrt{-1}/m')\in\mu_{m'}$ by the {\it transition morphism of projective system}
$$\mu_m\ni\la\,\mapsto\,\la^{m/m'}\in\mu_{m'}\q(m/m'\in\Z).$$
\par\nin The action of $T_s$ is then defined to be the action of the canonical generator for $m$ sufficiently divisible. So the above Grothendieck ring may be denoted also by $K_0^{T_s}({\rm Var}_S)^{\sim}$.
\par\htt{A.2}{}\msn
{\bf A.2.~Motivic descent theorem.} In this section we assume $f$ is a polynomial, and $U_0$ in \hl{2.1}{2.1} is a Zariski-open subset of $X$. In the notation of (\hl{2.1}{2.1}--\hl{2.2}{2}) and \hl{A.1}{A.1}, we have a canonical $T_s$-action on $\Xt_{\si}^{\circ}$ over $\X_{\si}^{\circ}$ ($\si\les\Gf$) as is explained after (\hl{2.2.2}{2.2.2}). Similarly we have the $T_s$ action on $\Yt_{\xi}^{\circ}$ over $\Y_{\xi}^{\circ}$ ($\xi\in\Xi'$) in the notation at the beginning of Section\,\,\hl{2.3}{2.3}. We have the trivial action of $T_s$ on $\X_{\si}^{\rm pr}$, $\Y_{\xi}^{\rm pr}$ for $\si\les\Gf$, $\xi\in\Xi'$.
Here $\Xt_{\si}^{\circ}$,$\X_{\si}^{\rm pr}$, etc.\ denote algebraic varieties (instead of analytic spaces) as is explained at the beginning of \hl{A.1}{A.1}.
\sk
Define the {\it motivic nearby fibers\1} of $f$ at $0$ in $K_0^{\muh}({\rm Var}_{\X_0^{\rm ex}})^{\sim}$, $K_0^{\muh}({\rm Var}_{\Y_0^{\rm ex}})^{\sim}$ respectively by
\htt{A.2.1}{}
$$\aligned\Sc_{f,\X_0^{\rm ex}}&:=\msum_{\si\les\Gf}\,(1{-}\Lb)^{k(\si)-d(\si)}\bl([\Xt_{\si}^{\circ}]+(1{-}\Lb)[\X_{\si}^{\rm pr}]\br),\\ \Sc_{f,\Y_0^{\rm ex}}&:=\msum_{\xi\in\Xi'}\,(1{-}\Lb)^{d_{\xi}-\kc_{\xi}-1}\bl([\Yt_{\xi}^{\circ}]+(1{-}\Lb)[\Y_{\xi}^{\rm pr}]\br),\endaligned
\leqno{\rm(A.2.1)}$$
\par\nin see Remark\,\,\hl{R2.2b}{2.2b} for the compatibility with Theorem\,\,\hl{T1}{1}.
By definition the image of $\Sc_{f,\Y_0^{\rm ex}}$ in $K_0^{\muh}({\rm Var}_{\C})^{\sim}$ via the direct image by $\pi\ssc\rho$ coincides with the one defined by Denef and Loeser \cite{DL}, \cite{Loe} in $K_0^{\muh}({\rm Var}_{\C})^{\sim}$, see a remark after (\hl{2.2.2}{2.2.2}). (The direct image is defined by composing the morphisms to $\Y_0^{\rm ex}$ with $\Y_0^{\rm ex}\to \{0\}$.)
\par\htt{TA.2}{}\msn
{\bf Theorem\,\,A.2.} {\it We have the equality
\htt{A.2.2}{}
$$\rho_!\1\Sc_{f,\Y_0^{\rm ex}}=\Sc_{f,\X_0^{\rm ex}}\q\h{in}\,\,\,\,K_0^{\muh}({\rm Var}_{\X_0^{\rm ex}})^{\sim},
\leqno{\rm(A.2.2)}$$
\par\nin where the direct image $\rho_!$ is defined by composing the morphisms to $\Y_0^{\rm ex}$ with $\rho:\Y_0^{\rm ex}\to\X_0^{\rm ex}$.}
\msn
{\it Proof.} For $\xi\in\Xi'$, there is $\si\les\Gf$ such that the corresponding $\eta_{\si}\in\Si'$ coincides with $\eta_{\xi}$ (the minimal dimensional cone containing $\xi$, see Remark\,\,\hl{R2.1d}{2.1d}), and the morphism $\rho$ induces the surjective morphism
\htt{A.2.3}{}
$$\rho_{\xi}:\Y_{\xi}\onto\X_{\si},
\leqno{\rm(A.2.3)}$$
\par\nin whose fibers are tori, see also (\hl{2.1.2}{2.1.2}--\hl{2.1.3}{3}). In view of Proposition\,\,\hl{P2.3}{2.3} and Corollary\,\,\hl{C2.3}{2.3}, the proof of Theorem\,\,\hl{TA.2}{A.2} is then reduced to the following.
\par\htt{PA.2}{}\msn
{\bf Proposition\,\,A.2.} {\it In the above notation, we have the equalities in} $K_0^{\muh}({\rm Var}_{\X_0^{\rm ex}})^{\sim}$\,:
\htt{A.2.4}{}
$$\aligned\rho_!\1[\Yt_{\xi}^{\circ}]&=(\Lb{-}1)^{d_{\eta_{\si}}-d_{\xi}}\1[\Xt_{\si}^{\circ}],\\ \rho_!\1[\Y_{\xi}^{\rm pr}]&=(\Lb{-}1)^{d_{\eta_{\si}}-d_{\xi}}\1[\X_{\si}^{\rm pr}].\endaligned
\leqno{\rm(A.2.4)}$$
\par\nin \msn
{\it Proof.} We proceed by induction on $d_{\rm rel}:=d_{\eta_{\si}}{-}\,d_{\xi}$. In the case $d_{\rm rel}=0$, the assertion follows from the compatibility of the covering spaces explained before (\hl{2.3.4}{2.3.4}), since (\hl{A.2.3}{A.2.3}) is an isomorphism by (\hl{2.1.2}{2.1.2}--\hl{2.1.3}{3}) when $d_{\rm rel}=0$.
\sk
In the case $d_{\rm rel}>0$, take $\ga\in\Xi'$ such that $\xi<\ga\subset\eta_{\si}$ and $d_{\ga}=d_{\xi}+1$. Let $\xi_1<\ga$ be the 1-dimensional face of $\ga$ not contained in $\xi$. Restricting to the half space $\xi^{\perp}\cap\xi_1^{\vee}\subset\xi^{\perp}$, we get the $\C$-subalgebra
$$R_{\xi,\ga}:=\C[\Z^n\cap\xi^{\perp}\cap\xi_1^{\vee}]\,\,\subset\,\,R_{\xi}:=\C[\Z^n\cap\xi^{\perp}],$$
\par\nin together with the surjection
$$R_{\xi,\ga}\,\onto\,R_{\ga}:=\C[\Z^n\cap\ga{}^{\perp}]=\C[\Z^n\cap\xi^{\perp}\cap\xi_1^{\perp}].$$
\par\nin These morphisms of $\C$-algebras correspond to the open and closed immersions of varieties
$$\Y_{\xi}\,\into\,\Y_{\xi,\ga}\,\hookleftarrow\,\Y_{\ga},$$
\par\nin with $\Y_{\xi,\ga}:={\rm Spec}\,R_{\xi,\ga}$. Set-theoretically we have $\Y_{\xi,\ga}=\Y_{\xi}\sqcup\Y_{\ga}$.
\sk
There is moreover a natural inclusion $R_{\ga}\into R_{\xi,\ga}$, which is equivalent to the smooth surjective morphism
$$\ph_{\ga,\1\xi}:\Y_{\xi,\ga}\,\onto\,\Y_{\ga},$$
\par\nin with fibers $\Ab^1$. This is a line bundle with zero-section given by the above closed immersion. We denote by $p_{\ga,\1\xi}:\Y_{\xi}\onto\Y_{\ga}$ the restriction of $\ph_{\ga,\1\xi}$ to $\Y_{\xi}$.
\sk
Since (\hl{A.2.4}{A.2.4}) holds for $\ga,\si$ by inductive assumption, the proofs of Theorem\,\,\hl{TA.2}{A.2} and Proposition\,\,\hl{PA.2}{A.2} are then reduced to the following.
\par\htt{LA.2}{}\msn
{\bf Lemma\,\,A.2.} {\it In the above notation, we have the equalities in} $K_0^{\muh}({\rm Var}_{\Y_{\ga}})^{\sim}$\,:
\htt{A.2.5}{}
$$\aligned(p_{\ga,\1\xi})_!\1[\Yt_{\xi}^{\circ}]&=(\Lb{-}1)\1[\Yt_{\ga}^{\circ}],\\(p_{\ga,\1\xi})_!\1[\Y_{\xi}^{\rm pr}]&=(\Lb{-}1)\1[\Y_{\ga}^{\rm pr}].\endaligned
\leqno{\rm(A.2.5)}$$
\par\nin \msn
{\it Proof.} The last equality of (\hl{A.2.5}{A.2.5}) follows from the relation (\hl{A.2.1}{A.2.1}), since the fibers of $p_{\ga,\1\xi}$ are 1-dimensional tori, and we have $\Y_{\xi}^{\rm pr}=p_{\ga,\1\xi}^{-1}(\Y_{\ga}^{\rm pr})\subset\Y_{\xi}$.
\sk
As for the first equality of (\hl{A.2.5}{A.2.5}), there is a (unique) normal variety $\Yt_{\xi,\ga}^{\circ}$ which is finite over $\Y_{\xi,\ga}^{\circ}:=\Y_{\xi,\ga}\setminus\Y_0^{\rm pr}$ and whose restriction over $\Y_{\xi}^{\circ}$ is isomorphic to $\Yt_{\xi}^{\circ}$. Indeed, its affine ring can be given by the {\it integral closure\1} of the affine ring of $\Y_{\xi,\ga}^{\circ}$ in the affine ring of $\Yt_{\xi}^{\circ}$. It is also possible to consider a smooth variety $W$ having a proper morphism $q$ to $\Y_{\xi,\ga}^{\circ}$ and whose restriction over $\Y_{\xi}^{\circ}$ is isomorphic to $\Yt_{\xi}^{\circ}$ (using a desingularization of some partial compactification of $\Yt_{\xi}^{\circ}$). We then get $\Yt_{\xi,\ga}^{\circ}={\rm Spec}\,q_*\OO_W$ (which is the Stein factorization of $q$, see \cite[III, Corollary 11.5]{Ha}).
\sk
Fix $y\in\Y_{\ga}^{\circ}$. Let $\Yt_{\xi,\ga,y}^{\circ}$ be the restriction of $\Yt_{\xi,\ga}^{\circ}$ over $\ph_{\ga,\1\xi}^{-1}(y)\,(\1\cong\C^)$. By definition, this is naturally isomorphic to the normalization $\Zt_y$ of the curve
$$\Zc_y:=\bl\{(z,u)\in\C^2\,\big|\,z^{\de_{\si,\xi}}=u^{\de_{\si,\xi_1}}\1h_{\si,\1\ga}(y)\br\}.$$
\par\nin with $z,u$ the coordinates of $\C$ and $\ph_{\ga,\1\xi}^{\,\,-1}(y)\,(\cong\C)$ respectively, see Remark\,\,\hl{RA.2a}{A.2a} below. Here $\de_{\si,\1\xi}\,$, $\de_{\si,\1\xi_1}$ are as in Remark\,\,\hl{R2.2e}{2.2e}.
\sk
We see that the number of irreducible components of $\Zc_y$ is equal to $\de_{\si,\1\ga}$, since
\htt{A.2.6}{}
$$\de_{\si,\1\ga}={\rm GCD}(\de_{\si,\1\xi},\de_{\si,\1\xi_1}),
\leqno{\rm(A.2.6)}$$
\par\nin as a consequence of Remark\,\,\hl{R2.2e}{2.2e}. (Notice that the action of $\muh$ is not necessarily diagonal, for instance, if $\de_{\si,\1\xi}/\de_{\si,\ga}=\de_{\si,\ga}$.)
\sk
Passing to the normalization, we have the decomposition
$$\Zt_y=\msqcup_i\,\Zt_{y,i}\q\h{with}\q\Zt_{y,i}\cong\C\q\bl(i\in[1,\de_{\si,\1\ga}]\br).$$
\par\nin The canonical morphism $\Zt_{y,i}\to\ph_{\ga,\1\xi}^{\,\,-1}(y)\,(\cong\C)$ is identified with the morphism
$$\C\,\ni\,v\,\mapsto\,v^{\de'_{\xi,\ga}}\,\in\,\C,$$
\par\nin where $\de'_{\xi,\ga}:=\de_{\si,\xi}/\de_{\si,\ga}$.
The identification $\ph_{\ga,\1\xi}^{\,\,-1}(y)\cong\C$ is unique up to constant multiple, since $\ph_{\ga,\1\xi}$ is a line bundle. We thus get a unique structure of $\C$-vector space on each $\Zc_{y,i}$.
\sk
We now show that there is a canonical morphism $\Yt_{\xi,\ga}^{\circ}\to\Yt_{\ga}^{\circ}$ whose fibers are 1-dimensional $\C$-vector spaces.
(This is closely related to Remark\,\,\hl{R2.1g}{2.1g} and also Lemma\,\,\hl{L2.4}{2.4}.)
Let $\Ga'_{\xi}$ be the cyclic subgroup of order $\de'_{\xi,\ga}$ of the cyclic covering transformation group $\Ga_{\xi}$ of $p_{\xi}^{\circ}$ which has order $\de_{\si,\,\xi}$. This defines a quotient covering space $\Yt_{\xi}^{\circ\prime}$ of order $\de_{\si,\ga}$ of $\Yt_{\xi}^{\circ}\onto\Y_{\xi}^{\circ}$. This quotient covering space is {\it constant\1} along the fibers of
$$p_{\ga,\1\xi}:\Y_{\xi}\onto\Y_{\ga},$$
\par\nin and can be extended to an unramified covering space $\Yt_{\xi,\ga}^{\circ\prime}$ over $\Y_{\xi,\ga}^{\circ}$. Here the restriction of $\Yt_{\xi,\ga}^{\circ\prime}$ over $\Y_{\ga}^{\circ}$ is naturally isomorphic to $\Yt_{\ga}^{\circ}$ by definition. Indeed, the covering spaces over $\Y_{\xi}^{\circ}$, $\Y_{\ga}^{\circ}$ are defined by dividing (the pull-backs of) $f_{\si}$ by appropriate monomials (depending also on the covering degrees) as is explained before (\hl{A.2.6}{A.2.6}), see also Remark\,\,\hl{RA.2a}{A.2a} below.
(Note also that the action of $\Ga'_{\xi}$ on $\Yt_{\xi}^{\circ}$ is identified with that of the covering transformation group of $\Zt_{y,i}\to\ph_{\ga,\1\xi}^{\,\,-1}(y)$ for any $i\in[1,\de_{\si,\ga}]$ by the definition of $\Zc_y$.)
\sk
By the above argument, $\Yt_{\xi,\ga}^{\circ}$ is a line bundle over $\Yt_{\ga}^{\circ}$. The algebraicity of this line bundle can be verified by using the base change by a cyclic finite \'etale morphism $\Y_{\ga}^{\circ\prime\prime}\to\Y_{\ga}^{\circ}$ associated with the $\de_{\si,\ga}$\1th root of unity of $h_{\si,\ga}(y)$, and decomposing the direct image of the associated locally free sheaf by the action of the cyclic covering transformation group.
\sk
Since the complement of the zero-section in $\Zt$ is isomorphic to $\Yt_{\xi}^{\circ}$ over $\Y_{\ga}$ via $p_{\ga,\1\xi}$ by definition, we then get the equalities in $K_0^{\muh}({\rm Var}_{\Y_{\ga}})^{\sim}$
\htt{A.2.7}{}
$$(p_{\ga,\1\xi})_!\1[\Yt_{\xi}^{\circ}]=[\Zt]-[\Yt_{\ga}^{\circ}]=(\Lb{-}1)[\Yt_{\ga}^{\circ}],
\leqno{\rm(A.2.7)}$$
\par\nin where the last equality follows from the relation (\hl{A.1.2}{A.1.2}). Thus the first equality of (\hl{A.2.5}{A.2.5}) is also proved.
This finishes the proofs of Lemma\,\,\hl{LA.2}{A.2}, Proposition\,\,\hl{PA.2}{A.2}, and Theorem\,\,\hl{TA.2}{A.2}.
\par\htt{RA.2a}{}\msn
{\bf Remark\,\,A.2a.} In the argument before (\hl{A.2.6}{A.2.6}), the two covering spaces
$$\Yt_{\xi}^{\circ}\onto\Y_{\xi}^{\circ},\q\Yt_{\ga}^{\circ}\onto\Y_{\ga}^{\circ}$$
\par\nin can be defined respectively by the equations
$$z^{\de_{\si,\xi}}=h_{\si,\1\xi}(u,y),\q z^{\de_{\si,\ga}}=h_{\si,\1\ga}(y),$$
\par\nin with $h_{\si,\1\xi}(u,y):=u^{\de_{\si,\xi_1}}\1h_{\si,\1\ga}(y)$. We can obtain $h_{\si,\1\xi}(u,y)$, $h_{\si,\1\ga}(y)$ by dividing (the pull-backs of) $f_{\si}$ by appropriate monomials in a similar way to the definition of $h_{\si}$ in the beginning of \hl{2.2}{2.2} (where the monomial is identified with a power of a coordinate of a line bundle). In the case of $h_{\si,\1\xi}(u,y)$ and $h_{\si,\1\ga}(y)$, these monomials are given respectively by the monomials
$$\mprod_{i=2}^{d_{\ga}}\,y_i^{\de_{\si,\xi_i}},\q\mprod_{i=1}^{d_{\ga}}\,y_i^{\de_{\si,\xi_i}},$$
\par\nin corresponding to the 1-dimensional cones $\xi_i<\xi$ ($i\in[2,d_{\ga}]$) and $\xi_i<\ga$ ($i\in[1,d_{\ga}]$) so that we get the equalities
$$h_{\si,\1\xi}(u,y)=f_{\si}/\mprod_{i=2}^{d_{\ga}}\,y_i^{\de_{\si,\xi_i}},\q h_{\si,\1\ga}(y)=f_{\si}/\mprod_{i=1}^{d_{\ga}}\,y_i^{\de_{\si,\xi_i}}.$$
\par\nin Here $y_i:=x^{\nu^{(i)}}$ with $\nu^{(i)}\in\Z^n$, and the $\nu^{(i)}$ are extended to free generators of $\Z^n$ satisfying the orthonormal relation
$$\langle\mu^{(j)},\nu^{(i)}\rangle=\de_{i,j},$$
\par\nin with $\mu^{(j)}\in\Z^n\cap\xi_j$ ($j\in[1,d_{\ga}]$), choosing an $n$-dimensional cone in $\Xi'$ containing $\ga$ (since $\Xi$ is a {\it smooth\1} subdivision). We then get (see (\hl{2.1.3}{2.1.3}))
$$\aligned h_{\si,\1\xi}(u,y)&\in\C[\Z^n\cap\xi^{\perp}]=\Ga(\Y_{\xi},\OO_{\Y_{\xi}}),\\ h_{\si,\1\ga}(y)&\in\C[\Z^n\cap\ga^{\perp}]=\Ga(\Y_{\ga},\OO_{\Y_{\ga}}).\endaligned$$
\par\nin \sk
The variables $z$, $u$ in the definition of $\Zc_y$ before (\hl{A.2.6}{A.2.6}) are identified respectively with
$$\mprod_{i=2}^{d_{\ga}}\,y_i^{-\de_{\si,\xi_i}/\de_{\si,\xi}},\q y_1,$$
\par\nin and the equation in the definition of $\Zc_y$ is essentially the pull-back of the equation $1=f_{\si}$, see also Remark\,\,\hl{R2.1f}{2.1f}.
\par\htt{RA.2b}{}\msn
{\bf Remark\,\,A.2b.} It seems to be noted in some paper quoted by \cite{Sta} that one cannot prove the motivic nearby fiber formula for Newton non-degenerate functions in the Grothendieck ring of complex algebraic varieties with good $\muh$-action in terms of the associated toric variety in a similar way to the normal crossing case by Denef and Loeser \cite{DL} (and some formula seems to be stated in the Grothendieck ring of mixed Hodge structures with good $\muh$-action), although its reason does not seem very clear.
\sk
It is also unclear if there is an algorithm to compute the Euler characteristic Hodge numbers of toric hypersurfaces in the ``non-prime" case, since the non-middle cohomology of a quasi-smooth compactification of a toric hypersurface via a quasi-smooth subdivision of the dual fan may have Hodge level greater than 0, see \cite{ic}.


\begin{thebibliography}{KKMS\,73}
\bibitem[BBD\,82]{BBD} Beilinson, A., Bernstein, J., Deligne, P., Faisceaux pervers, Ast\'erisque 100, Soc.\ Math.\ France, Paris, 1982.
\bibitem[Bi\,05]{Bi} Bittner, F., On motivic zeta functions and the motivic nearby fiber, Math. Z. 249 (2005), 63--83.
\bibitem[Br\,70]{Br} Brieskorn, E., Die Monodromie der isolierten Singularit\"aten von Hyperfl\"achen, Manuscripta Math. 2 (1970), 103--161.
\bibitem[BuSa\,05]{BS1} Budur, N., Saito, M., Multiplier ideals, $V$-filtration, and spectrum, J.\ Alg.\ Geom.\ 14 (2005), 269--282.
\bibitem[BuSa\,10]{BS2} Budur, N., Saito, M., Jumping coefficients and spectrum of a hyperplane arrangement, Math.\ Ann.\ 347 (2010), 545--579.
\bibitem[Da\,78]{Da1} Danilov, V.I., The geometry of toric varieties, Russian Math.\ Surveys 33 (1978), 97--154.
\bibitem[Da\,79]{Da2} Danilov, V.I., Newton polytopes and vanishing cohomology, Funct.\ Anal.\ Appl.\ 13 (1979), 103--115.
\bibitem[DGPS\,19]{Sing} Decker, W., Greuel, G.-M., Pfister, G., Sch\"onemann, H., {\sc Singular} 4-1-2 --- A computer algebra system for polynomial computations, available at http://www.singular.uni-kl.de (2019).
\bibitem[DeLo\,01]{DL} Denef, J., Loeser, F., Geometry on arc spaces of algebraic varieties, in European Congress of Mathematics, Vol.\ 1 (Barcelona, 2000), Progr.\ Math.\ 201, Birkha\"user, Basel, 2001, 327--348.
\bibitem[DiSa\,14]{DS} Dimca, A., Saito, M., Some remarks on limit mixed Hodge structures and spectrum, An.\ \c Stiin\c t.\ Univ.\ ``Ovidius'' Constan\c ta Ser. Mat. 22 (2014), 69--78.
\bibitem[EhLo\,82]{EL} Ehlers, F., Lo, K.-C., Minimal characteristic exponent of the Gauss-Manin connection of isolated singular point and Newton polytope, Math.\ Ann.\ 259 (1982), 431--141.
\bibitem[GrRe\,58]{GR} Grauert, H., Remmert, R., Komplexe R\"aume, Math.\ Ann.\ 136 (1958), 245--318.
\bibitem[Gr\,61]{Gr1} Grothendieck, A., El\'ements de g\'eom\'etrie alg\'ebrique III-1, Publ.\ Math.\ IHES 11, 1961.
\bibitem[Gr\,71]{Gr2} Grothendieck, A., Rev\^etements \'etales et groupe fondamental, S\'eminaire de G\'eom\'etrie Alg\'ebrique 1, Lect.\ notes in Math.\ 224, Springer, Berlin, 1971.
\bibitem[Ha\,77]{Ha} Hartshorne, R., Algebraic Geometry, Springer, Berlin, 1977.
\bibitem[JKSY\,22]{JKSY1} Jung, S.-J., Kim, I.-K., Saito, M., Yoon Y., Hodge ideals and spectrum of isolated hypersurface singularities, Ann.\ Inst.\ Fourier (Grenoble) 72 (2022), 465--510 (Proposition and formula numbers are changed by the publisher from arxiv:1904.02453).
\bibitem[JKSY\,23]{JKSY2} Jung, S.-J., Kim, I.-K., Saito, M., Yoon Y., Spectrum of non-degenerate functions with simplicial Newton polytopes, arxiv:1911.09465v4.
\bibitem[KKMS\,73]{KKMS} Kempf, G., Knudsen, F., Mumford, D., Saint-Donat, B., Toroidal Embeddings I, Lect.\ Notes Math.\ 339, Springer, Berlin, 1973.
\bibitem[Ko]{Ko} Kouchinirenko, A.G., Poly\`edres de Newton et nombres de Milnor, Inv.\ Math.\ 32 (1976), 1--31.
\bibitem[Loe\,05]{Loe} Loeser, F., Seattle lectures on motivic integration, in Algebraic geometry - Seattle 2005, Proc.\ Symp.\ Pure Math.\ Vol.\ 80, Part 2, Am.\ Math.\ Soc., 2009, pp.~745--784.
\bibitem[Loo\,02]{Loo} Looijenga, E., Motivic measures, S\'eminaire Bourbaki, exp.\ 874, Ast\'erisque 276 (2002), 267--297.
\bibitem[Ma\,75]{Ma} Malgrange, B., Le polyn\^ome de Bernstein d'une singularit\'e isol\'ee, Lect.\ Notes in Math.\ 459, Springer, Berlin, 1975, pp. 98--119.
\bibitem[MSS\,13]{MSS} Maxim, L., Saito, M., Sch\"urmann, J., Hirzebruch-Milnor classes of complete intersections, Adv.\ Math.\ 241 (2013), 220--245.
\bibitem[Od\,88]{Od} Oda, T., Convex bodies and algebraic geometry, Springer, Berlin, 1988.
\bibitem[Sa\,88a]{exp} Saito, M., Exponents and Newton polytopes of isolated hypersurface singularities, Math.\ Ann.\ 281 (1988), 411--417.
\bibitem[Sa\,88b]{mhp} Saito, M., Modules de Hodge polarisables, Publ.\ RIMS, Kyoto Univ.\ 24 (1988), 849--995.
\bibitem[Sa\,90]{mhm} Saito, M., Mixed Hodge modules, Publ.\ RIMS, Kyoto Univ.\ 26 (1990), 221--333.
\bibitem[Sa\,91]{ste} Saito, M., On Steenbrink's conjecture, Math.\ Ann.\ 289 (1991), 703--716.
\bibitem[Sa\,94]{mic} Saito, M., On microlocal $b$-function, Bull.\ Soc.\ Math.\ France 122 (1994), 163--184.
\bibitem[Sa\,07]{bfun} Saito, M., Multiplier ideals, $b$-function, and spectrum of a hypersurface singularity, Compos.\ Math.\ 143 (2007), 1050--1068.
\bibitem[Sa\,16]{hi} Saito, M., Hodge ideals and microlocal $V$-filtration, arxiv:1612.08667.
\bibitem[Sa\,18]{wt0} Saito, M., Weight zero part of the first cohomology of complex algebraic varieties, arxiv:1804.03632.
\bibitem[Sa\,20]{ic} Saito, M., Intersection complexes of toric varieties and mixed Hodge modules, arxiv:2006.04081.
\bibitem[ScSt\,85]{SS} Scherk, J., Steenbrink, J.H.M., On the mixed Hodge structure on the cohomology of the Milnor fibre, Math.\ Ann.\ 271 (1985), 641--665.
\bibitem[SSS\,91]{SSS} Schrauwen, R., Steenbrink, J.H.M., Stevens, J., Spectral pairs and the topology of curve singularities, Proc.\ Symp.\ pure Math.\ Vol.\ 53, Am.\ Math.\ Soc., 1991, pp.~305--328.
\bibitem[Se\,55]{Se} Serre, J.-P., Faisceaux alg\'ebriques coh\'erents, Ann.\ Math.\ 61 (1955), 197--278.
\bibitem[Sta\,17]{Sta} Stapledon, A., Formulas for monodromy, Res.\ Math.\ Sci.\ 4 (2017), Paper No.\ 8.
\bibitem[Ste\,77]{St1} Steenbrink, J.H.M., Mixed Hodge structure on the vanishing cohomology, in Real and complex singularities, Sijthoff and Noordhoff, Alphen aan den Rijn, 1977, pp. 525--563.
\bibitem[Ste\,89]{St2} Steenbrink, J.H.M., The spectrum of hypersurface singularities, Ast\'erisque 179-180 (1989), 163--184.
\bibitem[Va\,76]{Va1} Varchenko, A.N., Zeta-function of monodromy and Newton's diagram, Inv.\ Math.\ 37 (1976), 253--262.
\bibitem[Va\,82]{Va2} Varchenko, A.N., Asymptotic Hodge structure in the vanishing cohomology, Math.\ USSR-Izv.\ 18 (1982), 469--512.
\bibitem[Ve]{Ve} Verdier, J.-L., Dualit\'e dans la cohomologie des espaces localement compacts, S\'eminaire Bourbaki exp.\ 300, 1966, pp.~337--349.
\end{thebibliography}
\end{document}